\newif\ifpreprint

\pdfoutput=1

\preprinttrue

\ifpreprint

\documentclass[12pt]{article}
\usepackage{abstract} 

\usepackage{titlesec}
\def\sectionfont{\sffamily\Large\bfseries\boldmath}
\def\subsectionfont{\sffamily\large\bfseries\boldmath}
\def\paragraphfont{\sffamily\normalsize\bfseries\boldmath}
\titleformat*{\section}{\sectionfont}
\titleformat*{\subsection}{\subsectionfont}
\titleformat*{\subsubsection}{\paragraphfont}
\titleformat*{\paragraph}{\paragraphfont}
\titleformat*{\subparagraph}{\paragraphfont}
\usepackage[small,labelfont={bf,sf}]{caption}  

\usepackage{natbib}
\usepackage{amsmath,amssymb,amsthm,enumitem, mathtools}

\newtheorem{theorem}{Theorem}[section]  

\else

\RequirePackage{fix-cm}
\documentclass[smallcondensed,natbib]{svjour3}     
\smartqed  

\fi

\usepackage{graphicx}

\usepackage{multirow}
\usepackage{algorithm}
\usepackage[noend]{algpseudocode}

\def\basiceval#1{\the\numexpr#1\relax}

\usepackage[
tight-spacing=true,
round-mode=places,
table-number-alignment=left,
table-text-alignment=left,
]{siunitx}
\usepackage{csvsimple}	
\usepackage{booktabs}           
\usepackage{adjustbox} 
\newlength{\tablelength}
\usepackage{tikz}
\usepackage{pgfplots}
\usetikzlibrary{matrix,chains,decorations.pathreplacing,arrows}
\usetikzlibrary{positioning,calc}
\usetikzlibrary{pgfplots.groupplots}
\definecolor{main}{RGB}{163, 31, 52}
\definecolor{lineColor1}{RGB}{163, 31, 52}
\definecolor{lineColor2}{RGB}{40, 141, 193}
\definecolor{lineColor3}{RGB}{238, 119, 0}
\usetikzlibrary{calc}
\usetikzlibrary{shapes.misc}
\definecolor{splitColor}{RGB}{202, 200, 200}
\definecolor{classColor1}{RGB}{163, 31, 52}
\definecolor{classColor2}{RGB}{4, 30, 65}
\definecolor{classColor3}{RGB}{0, 105, 143}
\definecolor{classColor4}{RGB}{0, 108, 103}
\definecolor{classColor5}{RGB}{238, 119, 0}
\definecolor{classColor6}{RGB}{246, 178, 33}
\tikzset{%
  neuron/.style={
    circle,
    draw,
    minimum size=.5cm
  },
  missing/.style={
    draw=none,
    scale=4,
    text height=0.333cm,
    execute at begin node=\color{black}$\vdots$
  },
}
\tikzset{class1/.style={draw, rounded rectangle, text=white, fill=classColor1!80}}
\tikzset{class2/.style={draw, rounded rectangle, text=white, fill=classColor2!80}}
\tikzset{class3/.style={draw, rounded rectangle, text=white, fill=classColor3!80}}
\tikzset{class4/.style={draw, rounded rectangle, text=white, fill=classColor4}}
\tikzset{class5/.style={draw, rounded rectangle, text=white, fill=classColor5}}
\tikzset{class6/.style={draw, rounded rectangle, text=white, fill=classColor6}}
\tikzset{split/.style={draw,rectangle,fill=splitColor!20}}
\tikzset{level 1/.style={sibling distance=14em},
         level 2/.style={sibling distance=8em},
         level 3/.style={sibling distance=3em}
        }
\tikzset{%
  every neuron/.style={
    circle,
    draw,
    minimum size=.8cm
  },
  neuron missing/.style={
    draw=none,
    scale=3,
    text height=0.333cm,
    execute at begin node=\color{black}$\vdots$
  },
}
\usetikzlibrary{arrows}
\tikzstyle{block}=[draw, align=center, inner sep=.5em, fill=main, text=white, minimum size=2em]
\tikzstyle{init} = [pin edge={to-,thin,black}]

\usepackage{hyperref}

\usepackage[xindy, acronyms]{glossaries}   

\definecolor{borderBoxTodo}{rgb}{0,0.309,0.6}
\definecolor{backgroundBoxTodo}{rgb}{0.701,0.854,1}
\definecolor{lineColorTodo}{rgb}{0.701,0.854,1}
\RequirePackage[colorinlistoftodos, backgroundcolor=backgroundBoxTodo, bordercolor = borderBoxTodo, linecolor = lineColorTodo]{todonotes}

\newcommand{\ie}{{\it i.e.}}

\newcommand{\ones}{\mathbf 1}
\newcommand{\reals}{{\mbox{\bf R}}}
\newcommand{\integers}{{\mbox{\bf Z}}}
\newcommand{\symm}{{\mbox{\bf S}}}  

\newcommand{\prob}{{\mathbf P}}
\newcommand{\distrib}{{\mathcal P}}

\newcommand{\tpose}{T}

\newcommand{\eps}{\epsilon}

\newcommand{\argmax}{\mathop{\rm argmax}}

\newcommand{\tree}{\mathcal{T}}

\newcommand{\strategy}{\mathcal{S}}
\newcommand{\intvars}{\mathcal{I}}
\newcommand{\tightconstraints}{\mathcal{T}}

\newcommand{\loss}{\mathcal{L}}

\newcommand{\Sec}{Section\;}

\newacronym{LO}{LO}{linear optimization problem}
\newacronym{QO}{QO}{quadratic optimization problem}
\newacronym{MIQO}{MIQO}{mixed-integer quadratic optimization problem}
\newacronym{MIO}{MIO}{mixed-integer optimization}
\newacronym{MILO}{MILO}{mixed-integer linear optimization problem}
\newacronym{MINLO}{MINLO}{mixed-integer nonlinear optimization problem}
\newacronym{SPM}{SPM}{successive projection method}
\newacronym{sBB}{sBB}{spacial branch and bound}
\newacronym{NLO}{NLO}{nonlinear optimization problem}
\newacronym{PWA}{PWA}{piecewise affine}
\newacronym{SVM}{SVM}{support vector machines}
\newacronym{OCT}{OCT}{optimal classification tree}
\newacronym[firstplural=optimal classification trees with-hyperplanes (OCT-Hs)]{OCT-H}{OCT-H}{optimal classification tree with-hyperplanes}
\newacronym{CART}{CART}{classification and regression tree}
\newacronym{NN}{NN}{neural network}
\newacronym{ReLU}{ReLU}{rectified linear unit}
\newacronym{MPC}{MPC}{model predictive control}
\newacronym{BNB}{B\&B}{branch-and-bound}
\newacronym{CP}{CP}{constraint program}

\begin{document}

\title{\bfseries \sffamily The Voice of Optimization}

\ifpreprint
\author{Dimitris Bertsimas and Bartolomeo Stellato}
\fi

\ifpreprint

\maketitle

\begin{abstract}
  We introduce the idea that using \glspl{OCT} and \glspl{OCT-H}, interpretable machine learning algorithms developed by~\citet{bertsimas2017a,bd_book}, we are able to obtain insight on the strategy behind the optimal solution in continuous and mixed-integer convex optimization problem as a function of key parameters 
that  affect the problem.
In this way,  optimization is not a black box anymore.
Instead, we redefine optimization as a multiclass classification problem where the predictor gives insights on the logic behind the optimal solution.
In other words, \glspl{OCT} and \glspl{OCT-H} give optimization a voice.
We show on several realistic examples that the accuracy behind our method is in the 90\%-100\% range, while even when the predictions are not correct, the degree of suboptimality or infeasibility is very low.
We compare optimal strategy predictions of  \glspl{OCT} and \glspl{OCT-H} and feedforward \glspl{NN} 
and conclude that the performance of \glspl{OCT-H}
 and NNs is comparable.  \glspl{OCT} are somewhat weaker but often  competitive.
Therefore, our approach provides a novel insightful understanding of optimal strategies to solve a broad class of continuous and mixed-integer optimization problems.


\end{abstract}

\fi

\glsresetall

\section{Introduction}

Optimization has a long and distinguished history that has had and continues to have genuine
impact in the world. In a typical optimization problem in the real world, practitioners see optimization as a black-box tool where they formulate the problem and they pass it to a
solver to find an optimal solution. Especially in high dimensional problems typically encountered in real world applications,
it is currently   not possible to interpret  or intuitively understand the optimal solution.
However,  independent from the  optimization algorithm  used, practitioners  would like to understand how problem  parameters affect the optimal decisions
in order  to get  intuition and interpretability behind the optimal solution.
Moreover,  in almost all   real world applications of optimization  the main objective is not to solve just one problem but to solve multiple similar instances
 that vary slightly from each other. In fact, in most real-world applications we solve similar optimization problems multiple times with varying data depending on problem-specific parameters.

Our goal  in this paper is to propose a framework to predict the optimal solution as parameters of the problem vary
and do so in an interpretable way.
A naive approach could be to learn directly the optimal solution from the problem parameters.
However, this method is computationally {\em intractable} and {\em imprecise}:
intractable, because it would require a potentially high dimensional predictor with hundreds of thousands of components depending on the decision variable size;
imprecise, because it would involve a regression task that would naturally carry a generalization error leading to suboptimal or even infeasible solutions.
Instead, our approach encodes the optimal solution with a small amount of information that we denote as {\em strategy}.
Depending on the problem class, a strategy can have different meanings but it always corresponds to the complete information needed to efficiently recover the optimal solution.

In recent times,  machine learning has also had significant impact to the world.
Optimization methods has been a major driver of its success~\citep{hastie2009,bd_book}.
In this paper, we apply machine learning to optimization with the objective to give a voice to optimization,
that is to provide interpretability and intuition behind optimal solutions.

First, we solve an optimization problem  for many parameter  combinations and obtain the optimal strategy: the set of active constraints as well as the value of the discrete variables.
Second, we encode the strategies with unique integer scalars. In this way, we have a mapping from parameters to the optimal solution strategies, which gives rise to a multiclass classification problem~\citep{hastie2009}.
We solve the multiclass classification problem using~\glspl{OCT} and \glspl{OCT-H} developed by \citet{bertsimas2017a,bd_book}, interpretable state of the art classification algorithms as well as \glspl{NN}, which, while not interpretable, serve as a benchmark to compare the accuracy of \glspl{OCT} and \glspl{OCT-H}.

\subsection{Related Work}%
\label{sub:related_work}

There has been a significant interest from both the computer science and the optimization communities to systematically analyze and solve optimization problems using machine learning.
From the first works applying machine learning as a substitute for mathematical optimization algorithms in the 1990s~\citep{smith1999}, this research direction has been increasingly active until recent results on learning for combinatorial optimization~\citep{bengio2018}.

\paragraph{Learning to tune algorithms.}
Most optimization methods consist of iterative routines with several parameters obtained through experts' knowledge and manual tuning.
For example, Gurobi optimizer has 157 parameters~\citep{gurobioptimization2020} that are carefully tuned in every release.
However, this tuning can be, on the one hand, very complex because it concerns several aspects of the problem that are not known a priori and, on the other hand, suboptimal for the instances considered in a specific application.
To overcome these limitations, the machine learning community studied the \emph{algorithm configuration} problem, \ie, the problem of automatically tuning algorithms from the problem instances of interest.
One of the first efficient automatic configuration methods is ParamILS~\citep{hutter2009} which iteratively improves the algorithm parameters using local search methods.
\cite{hutter2011} extended this idea in SMAC, a framework to learn the relationship between algorithm configurations and performance.
The proposed scheme relies on Gaussian Processes and Random Forests to perform the predictions.
Later, \cite{lopez2016} introduced irace, an algorithm that iteratively refines the candidate parameters configurations to identify the best ones for a specific application.
Despite their relevance in practice, these approaches aim to identify or predict the best algorithm parameters, independently from the algorithm ultimate task, \ie, they do not only consider optimization algorithms.
In addition, they do not consider interpretability of the predictors. For instance, SMAC could benefit from interpretable predictors such as~\glspl{OCT} to allow the user to understand why some parameter configurations are better than others.
Our work, instead, studies the relationship between the problem instances and the optimal solution focusing specifically on optimization algorithms and the interpretability of the optimal solutions.

%
%
%
%
%

\paragraph{Learning heuristics.}
Optimization methods not only rely on careful parameter tuning, but also on efficient heuristics.
For example, \gls{BNB} involves several heuristic decisions about the branching behavior that are hand-tuned into the solvers. However, heuristics tuning can be very complex because it concerns several aspects of the problem that are not known a priori.
\cite{khalil2016} propose to learn the branching rules showing performance improvements over commercial hand-tuned algorithms.
Similarly,~\cite{alvarez2017} approximate strong branching rules with learning methods.
Machine learning has been useful also to select reformulations and decompositions for~\gls{MIO}.
\citet{bonami2018} learn in which cases it is more efficient to solve~\gls{MIQO} by linearizing or not the cost function. They model it as a classification problem showing advantages compared to how this choice is made heuristically inside state-of-the-art solvers.
\citet{kruber2017} propose a similar method applied to decomposition selection for~\gls{MIO}.

\paragraph{Reinforcement learning for optimization.}
Another interesting line of research models optimization problems as control tasks to tackle using reinforcement learning~\citep{sutton2018}. Problems suitable to this framework include knapsack-like or network problems with multistage decisions.
\citet{dai2017} develop a method to learn heuristics over graph problems. In this way the node selection criterion becomes the output of a specialized neural network that does not depend on the graph size~\citep{dai2016}.
Every time a new node is visited,~\citet{dai2017} feed a graph representation of the problem to the~\gls{NN} obtaining a criterion suggesting the next node to select in the optimal solution.

\paragraph{Learning constraint programs.}
Constraint programming is a paradigm to model and solve combinatorial optimization problems very popular in the computer science community.
The first works on applying machine learning to automatically configure~\glspl{CP} date back to the 1990s~\citep{minton1996}.
Later,~\cite{clarke2002} used Decision Trees to replace computationally hard parts of counterexample guided SAT solving algorithms.
More recently,~\cite{xu2008} describe SATzilla, an automated approach for learning which candidate solvers are best on a given instance.
SATzilla won several SAT solver competitions because if its ability to adapt and pick the best algorithm for a problem instance.
Given a SAT instance, instead of solving it with different algorithms, SATzilla relies on a \emph{empirical hardness model} to predict how long each algorithm should take.
This model consists of a ridge regressor~\citep{hastie2009} after nonlinear transformation of the problem features.
\cite{selsam2018} applied recent advances in \gls{NN} archtectures to \glspl{CP} by directly predicting the solution or infeasibility.
Even though this approach did not give as good results as state-of-the-art methods, it introduces a new research direction for solving \glspl{CP}.
Therefore, the constraint programming community is also working on data-driven methods to improve the performance and understanding of solution algorithms.

%
%
%
%
%

\paragraph{Learning parametric programs.}
Even though recent approaches for integrating machine learning and optimization show promising results, they do not consider the parametric nature of the problems appearing in real-world applications.
It is often the case that practitioners solve the same problem with slightly varying parameters multiple times generating a large amount of data describing how the parameters affect the optimal solution.
There are only a few recent papers exploiting this information to build better solution algorithms.
The idea of learning the set of active constraints for parametric online optimization has been proposed by~\citet{misra2019}.
The authors frame the learning procedure as a sampling scheme where they collect all the active sets appearing from the parameters.
However, we found several limitations of  their approach.
First, in the online phase, they evaluate all the relevant active sets in parallel and pick the best one~\citep[\emph{ensemble policy}]{misra2019}.
Therefore, they do not solve the optimization problem as a multiclass classification problem and they are not able to gain insights on how the parameters affect the optimal strategy.
In addition, they do not tackle mixed-integer optimization problems but only continuous convex ones.
Finally, the sampling strategy by~\citet{misra2019} has to be tuned for each specific problem since it depends on at least four different parameters. This is because the authors compute the probabilistic guarantees based on how many new active sets appear over the samples in a window of a specific size~\citep[Section 3]{misra2019}.
In this work, instead, we provide a concise Good-Turing estimator~\citep{good1953} for the probability of finding new unseen strategies which can be directly applied to many different problem instances.
In the field of~\gls{MPC},~\citet{klauco2019} warm-start an online active set method for solving~\glspl{QO}.
However, in that work there is no rigorous sampling scheme with probability bounds to obtain the different active sets.
Instead, the authors either simulate the dynamical controlled system or provide an alternative gridding method to search for the relevant active sets.
Furthermore, the method by~\citet{klauco2019} is tailored to a specific linear control problem in the form of \gls{QO} and cannot tackle general convex or mixed-integer convex problems.
Learning for parametric programs can also speedup the online solution algorithms.
\citet{bertsimas2019} apply the framework in this paper to online mixed-integer optimization. By focusing on speed instead of interpretability, they obtain significant speedups compared to state-of-the-art solvers.

\paragraph{Sensitivity analysis.}
The study of how changes in the problem parameters affect the optimal solution has for long been studied in sensitivity analysis,~\citep[Chapter 5]{bertsimas1997} and~\citep[Section 5.6]{boyd2004} for introductions on the topic.
While sensitivity analysis is related to our work as it analyzes the effects of changes in problem parameters, it is fundamentally different both in philosophy and applicability.
In sensitivity analysis the problem parameters are uncertain and the goal is to understand how their perturbations affect the optimal solution. This aspect is important when, for example, the problem parameters are not known with high accuracy and we would like to understand how the solution would change in case of perturbations.
In this work instead, we consider problems without uncertainty and use previous data to learn how the problem parameters affect the optimal solution. Therefore, our problems are deterministic and we are not considering perturbations around any nominal value. As a matter of fact, the data we use for training are not restricted to lie close to any nominal point.
In addion, sensitivity analysis usually studies continuous optimization problems since it relies on the dual variables at the optimal solution to determine the effect of parameter perturbations. This is why there has been only limited work on sensitivity analysis for~\gls{MIO}. In contrast, we show that our method can be directly applied to problems with integer variables.

\subsection{Contributions}%
\label{sub:contributions}

In this paper, we propose a learning framework to give a voice to continuous and mixed-integer convex optimization problems.
With our approach we can reliably sample the occurring strategies using the Good-Turing estimator, learn an interpretable classifier using \glspl{OCT} and \glspl{OCT-H}, interpret the dependency between the optimal strategies and the key parameters from the resulting tree and solve the optimization problem using the learned predictor.

Our specific contributions include:
\begin{enumerate}
\item We introduce a new framework for gaining insights on the solution of optimization problems as a function of their key parameters. The optimizer becomes an interpretable machine learning classifier using \glspl{OCT} and \glspl{OCT-H} which highlights the relationship between the optimal solution strategy and the problem parameters. In this way optimization is no longer a black box and our approach gives it a voice that provides intuition and interpretability on the optimal solution.
\item We show that our method can be applied to a broad collection of optimization problems including convex continuous and convex mixed-integer optimization problem.
We do not pose any assumption on the dependency of the cost and constraints on the problem parameters.
\item We introduce a new exploration scheme based on the Good-Turing estimator~\citep{good1953} to discover the strategies arising from the problem parameters. This scheme allows us to reliably bound the probability of encountering unseen strategies in order to be sure our classifier accurately represents the optimization problem.
\item
In several realistic examples we show that the sample accuracy of our method is in the 90\%-100\% range, while even in the cases where the prediction is not correct the degree of suboptimality or infeasibility is very low.
We also compare the performance of \glspl{OCT} and \glspl{OCT-H} to \glspl{NN} which can achieve state-of-the-art performance across a wide variety of prediction tasks.
In our experiments we obtained comparable out-of-sample accuracy with \glspl{OCT}, \glspl{OCT-H} and \glspl{NN}.
\end{enumerate}
In other words, our approach provides a novel, reliable and insightful framework for understanding the strategies to compute the optimal solution of a broad class of continuous and mixed-integer optimization problems.

\subsection{Paper Structure}%
\label{sub:structure}
The structure of the  paper is as follows.
In Section \ref{sec:voice}, we define   an  optimal strategy to solve continuous and mixed-integer optimization problems and present several concrete examples   that demonstrate what we call the voice of optimization.
In Section \ref{sec:machine_learning}, we outline  the core part of our approach of using multiclass classification to learn the   mapping from parameters  to optimal strategies.
We further present an approach   to estimate how likely it is to encounter a parameter that leads to an optimal strategy that has not yet been seen.
In Section \ref{sec:machine_learning_optimizer}, we outline  our
Python  implementation  MLOPT (Machine Learning Optimizer).
In Section
\ref{sec:benchmarks}, we test our  approach on multiple examples from continuous and mixed-integer optimization and present the accuracy of predicting the optimal strategy of \glspl{OCT} and \glspl{OCT-H} in comparison with a \glspl{NN} implementation.
Section \ref{sec:conclusions} summarizes our conclusions.
Appendices~\ref{sec:trees} and~\ref{sec:neural_networks} briefly present  optimal classification trees (\glspl{OCT}  and \glspl{OCT-H})  and \glspl{NN}, respectively to make the paper self-contained.




\section{The Voice of Optimization}
\label{sec:voice}

In this section, we introduce the notion of  an   optimal strategy to solve continuous and mixed-integer optimization problems.

Given a parametric optimization problem, we define \emph{strategy}~$s(\theta)$ as the complete information needed to efficiently compute its optimal solution given the parameter $\theta \in \reals^p$.
We assume that the problem  is  always feasible  for every encountered value of~$\theta$.

Note that in the illustrative examples in this section we omit the details of the learning algorithm  which we explain in Section~\ref{sec:machine_learning}.
Instead, we focus on the interpretation of the resulting strategies which correspond to a decision tree.
For these examples the  accuracy of our approach to find the optimal solution  was always $100\%$, confirming that the resulting models accurately recover the optimal solution.
For simplicity of exposition, we sample the problem parameters in the examples of this section in simple regions such as intervals or balls around specific points but this is not a strict requirement for our approach as it will be more clear in the next sections.
In the examples in this section we used~\glspl{OCT} since their accuracy was 100\%  and they are more interpretable than~\glspl{OCT-H}.
In the last example, we used both an OCT and an OCT-H for comparison.

\subsection{Optimal Strategies in  Continuous Optimization}
\label{sub:continuous}

Consider the continuous optimization problem
\begin{equation}\label{eq:continuous_problem}
\begin{array}{ll}
\text{minimize} & f(\theta, x)\\
\text{subject to} & g(\theta, x) \le 0,\\
\end{array}
\end{equation}
where $x\in \reals^{n}$ is the vector of decision variables and $\theta \in \reals^p$ the vector of parameters affecting the problem data.
Functions $f:\reals^p \times \reals^n \to \reals$ and  $g: \reals^p \times \reals^n \to \reals^m$ are assumed to be convex in $x$.
Given a parameter $\theta$ we denote the optimal primal solution as $x^\star(\theta)$ and the optimal cost function value as $f(\theta, x^\star(\theta))$.

\paragraph{Tight constraints.}

Let us define the {\em tight constraints} $\tightconstraints(\theta)$ as the set of constraints that are satisfied as equalities at optimality,
\begin{equation}\label{eq:tight_constraints}
\tightconstraints(\theta) = \{i \in \{1, \dots, m\} \mid g_i(\theta, x^\star) = 0\}.
\end{equation}
Given the tight constraints, all the other constraints are no longer needed to solve the original problem.

For non-degenerate problems, the tight constraints correspond to the {\em support constraints}, \ie, the set of constraints that, if removed, would allow a decrease in $f(\theta, x^\star(\theta))$~\citep[Definition 2.1]{calafiore2010}.
In the case of \glspl{LO} the support constraints are the linearly independent constraints defining a basic feasible solution~\citep[Definition 2.9]{bertsimas1997}.
An important property of support constraints is that they cannot be more than the dimension $n$ of the decision variable~\citep[Proposition 1]{hoffman1979},~\citep[Lemma 2.2]{calafiore2010}.
This fact plays a key role in our method to reduce the complexity to predict the solution of parametric optimization problems.
The benefits are more evident when the number of constraints is much larger than the number of variables, \ie, $n \ll m$.

\paragraph{Multiple optimal solutions and degeneracy.}
In practice we can encounter problems with multiple optimal solutions or degeneracy.
When the cost function is not strongly convex in $x$, we can have multiple optimal solutions. In these cases we consider only one of the optimal solutions to be $x^\star$ since it is enough for our training purposes. Note that most solvers such as~\citep{gurobioptimization2020} return anyway only one solution and not the complete set of optimal solutions.
With degenerate problems, we can have more tight constraints than support constraints for an optimal solution~$x^\star$. This is because the support constraints set is no longer unique in case of degeneracy. However, we use the set of tight constraints which remains unique since it includes by definition all the constraints that are satisfied as equalities independently from being support constraints or not. In addition, our method is still efficient because the number of tight constraints is in practice much lower than the total number of constraints, even in case of degeneracy.
Therefore, we can directly apply our framework to problems with multiple optimal solutions and affected by degeneracy.

\paragraph{Solution strategy.}
We can now define our strategy as the index of tight constraints at the optimal solution, \ie, $s(\theta) = \tightconstraints(\theta)$.
Given the optimal strategy, solving~\eqref{eq:continuous_problem} corresponds to solving
\begin{equation}\label{eq:continuous_reduced_problem}
\begin{array}{ll}
\text{minimize} & f(\theta, x)\\
        \text{subject to} & g_i(\theta, x) \le 0, \quad \forall i \in \tightconstraints(\theta).
\end{array}
\end{equation}
This problem is easier than~\eqref{eq:continuous_problem}, especially when $n \ll m$.
Note that we can enforce the components $g_i$ that are linear in $x$ as equalities.
This further simplifies~\eqref{eq:continuous_reduced_problem} while preserving its convexity.
In case of~\gls{LO} and~\gls{QO} when the cost $f$ is linear or quadratic and the constraints $g$ are all linear, the solution of \eqref{eq:continuous_reduced_problem} corresponds to solving a linear system of equations defined by the KKT conditions~\cite[Section 10.2]{boyd2004}.

\subsubsection*{Inventory management}%
\label{ssub:inventory_example}

Consider an inventory management problem with horizon ${t = 0,\dots, T-1}$ with $T=30$. The decision variables are $u_t$, describing how much we order at time $t$ and $x_t$, describing the inventory level at time $t$.
The cost of ordering is $c=2$, $h=1$ is the holding cost and $p=3$ is the shortage cost.
We define the maximum quantity we can order each time as $M=3$.
The parameters are the product demand $d_t$ at time~$t$ and the initial value of the inventory~$x_{\rm init}$.
The optimization problem can be written as follows
\begin{equation}\label{eq:inventory_continuous_example}
\begin{array}{ll}
\text{minimize} & \displaystyle  \sum_{t=0}^{T-1} \max(h x_t, -p x_t) + c u_t\\
        \text{subject to} & x_{t+1} = x_{t} + u_{t} - d_{t}, \quad t=0,\dots, T-1\\
         & x_0 = x_{\rm init}\\
         & 0 \le u_t \le M, \quad t=0,\dots, T-1. \\
\end{array}
\end{equation}
Depending on $d_t$ with $t=0,\dots, T-1$ and $x_{\rm init}$, we need to adapt our ordering policy to minimize the cost.
We assume that  $d_t\in [1,3]$ and   $x_{\rm init}\in[7,13]$.

\begin{figure}
    \centering
    \begin{tikzpicture}[sibling distance=8em]
    \node[split] {$x_{\rm init} < 9.97$}
        child { node[split] {$x_{\rm init} < 7.91$}
                    child {node[class4]{4}}
                    child {node[class2]{2}}
                edge from parent node[left, yshift=1ex] {True}}
        child { node[split] {$d_5 < 1.91$}
                    child { node[split] {$x_{\rm init} < 11.61$}
                        child {node[class1]{1}}
                        child {node[class3]{3}}}
                    child { node[split] {$x_{\rm init} < 12.2$}
                        child {node[class1]{1}}
                        child {node[class3]{3}}}
                edge from parent node[right, yshift=1ex] {False}}
    ;
\end{tikzpicture}
    \caption{Decision strategy for the inventory example.
    \emph{Strategy 1}: do not order for the first 4 time steps, then order matching the demand. \emph{Strategy 2}: do not order for the first 5 time steps, then order matching the demand. \emph{Strategy 3}: do not order for the first 6 time steps, then order matching the demand. \emph{Strategy 4}: do not order for the first 3 time steps, then order matching the demand.}
    \label{fig:example_inventory_continuous}
\end{figure}
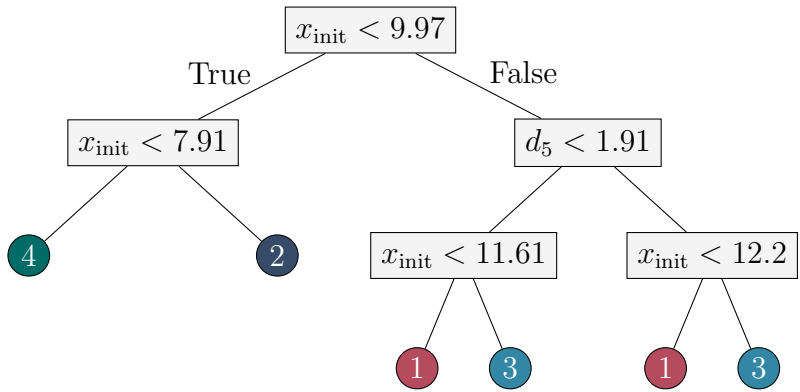

The strategy selection is summarized in Figure~\ref{fig:example_inventory_continuous} and can be easily interpreted: $u_t = 0$ for $t\leq t_0$  and then $u_t = d_t$ for $t> t_0$.
We can explain this rule from the problem formulation. As discussed before, the strategy tells us which constraints are tight, in this case when $u_t = 0$ and, therefore, when we do not order. In the other time steps, we can ignore the inequality constraints and our goal is to minimize the cost by matching the demand. In this way, we do not have to store anything in the inventory, \ie, $x_t = 0$ and $\max(hx_t, -px_t) = 0$. Note that we anyway need to pay the cost of ordering $cu_t$ since we have the satisfy the demand over the horizon.

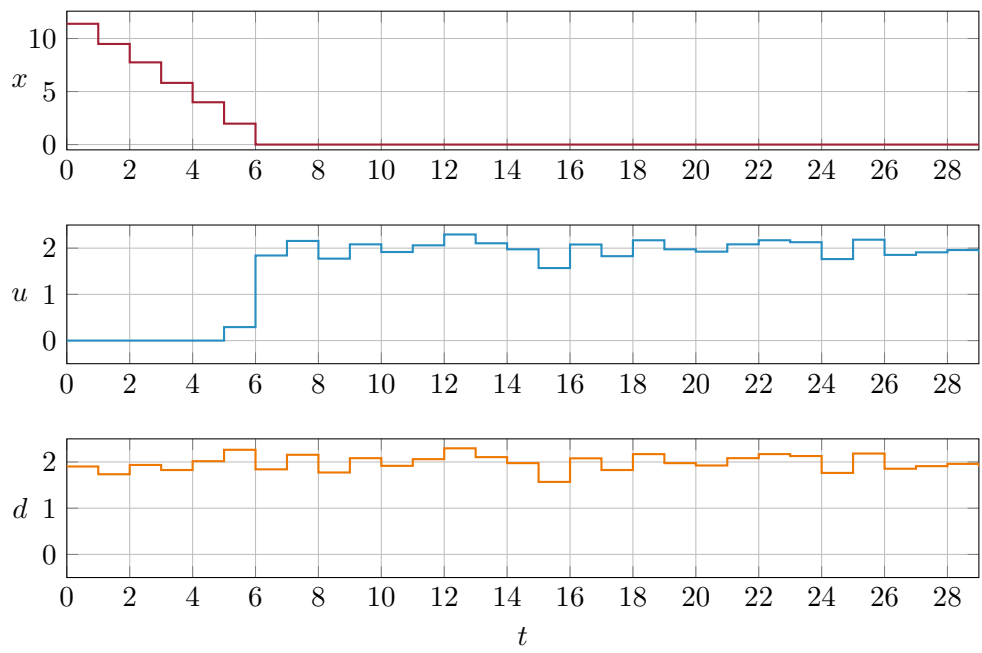
\begin{figure}
    \begin{center}
		\small
		\begin{tikzpicture}
        \begin{groupplot}[
        group style = {
        group size = 1 by 3,
        ylabels at = edge left
        },
    		    width=\columnwidth,
                height=0.25\columnwidth,
                grid=both, 
                xmin = 0.0,
                xmax = 29.,
                y label style={at={(axis description cs:0.05,.5)},rotate=-90},
				]
                \nextgroupplot[ylabel=$x$, ymin = -0.5]
				\addplot [lineColor1, const plot, thick]   table[x=t, y=x, col sep=comma] {./data/discussion/inventory/inventoryplot.csv};
                \nextgroupplot[ylabel=$u$, ymax = 2.5, ymin = -0.5]
				\addplot [lineColor2, const plot, thick]   table[x=t, y=u, col sep=comma] {./data/discussion/inventory/inventoryplot.csv};
                \nextgroupplot[xlabel={$t$}, ylabel=$d$, ymax = 2.5, ymin = -0.5]
				\addplot [lineColor3, const plot, thick]   table[x=t, y=d, col sep=comma] {./data/discussion/inventory/inventoryplot.csv};
            \end{groupplot}
		\end{tikzpicture}
    \end{center}
    \caption{Inventory behavior with Strategy 2. The lower bound on $u_t$ is active for the first 5 steps. Then $u_t = d_t$.}
    \label{fig:inventory_behavior}
\end{figure}

An inventory level trajectory example appears in Figure~\ref{fig:inventory_behavior}.
Let us outline the strategies depicted in Figure~\ref{fig:example_inventory_continuous}.
For example, if the initial inventory level $x_{\rm init}$ is below $7.91$, we should apply Strategy $4$, where we wait only $t_0=3$ time steps before ordering.
Otherwise, if $7.91 \le x_{\rm init} < 9.97$ we should wait for $t_0=5$ time steps before ordering because the initial inventory level is higher.
The other branches can be similarity interpreted.
Note that for this problem instance the decision is largely independent from $d_t$.
The only exception comes with $d_5$ that determines the choice of strategies in the right-hand side of the tree.

Note that this is a simple illustrative example and the strategies shown are not all the theoretical ones.
By allowing the problem parameters to take all the possible values we should expect a much larger number of strategies, and therefore a deeper tree.
However, even though real-world examples can be more complicated, we very rarely hit all the possible strategy combinations.

The strategies we derived  are related to the $(s, S)$ policies for inventory management~\citep{zheng1991}.
If the inventory level $x_t$ falls below value $s$, an $(s, S)$ policy orders a quantity $u_t$ so that the inventory is replenished to level $S$.
This policy is remarkably simple because it consists only in a {\em if-then-else} rule. The optimal decisions from our approach are as simple as the $(s, S)$ policies because we can describe them with the \emph{if-then-else} rules as in Figure~\ref{fig:example_inventory_continuous}.
However, our method gives better performance because it takes into account future predictions and can deal with more complex problem formulations with harder constraints.

\subsection{Optimal Strategies in Mixed-Integer Optimization}
\label{sub:mixed-integer}

When dealing with integer variables we address the following  problem
\begin{equation}\label{eq:integer_problem}
\begin{array}{ll}
\text{minimize} & f(\theta, x)\\
\text{subject to} & g(\theta, x) \le 0\\
&  x_\intvars \in \integers^d.
\end{array}
\end{equation}
where $\intvars$~is the set of indices for the variables constrained to be integer and $|\intvars| = d$.

\paragraph{Tight constraints.}
In this case the set of tight constraints does not uniquely define the optimal solution because of the integrality of some components of $x$.
However, when we fix the integer variables to their optimal values $x_{\intvars}^\star(\theta)$,~\eqref{eq:integer_problem} becomes a continuous convex optimization problem of the form
\begin{equation}\label{eq:integer_problem_continuous}
\begin{array}{ll}
\text{minimize} & f(\theta, x)\\
\text{subject to} & g(\theta, x) \le 0\\
&  x_\intvars = x_{\intvars}^\star(\theta).
\end{array}
\end{equation}
Note that the optimal cost of~\eqref{eq:integer_problem_continuous} and~\eqref{eq:integer_problem} are the same,
however,  optimal solutions may be different as there may be alternative optima.
After fixing the integer variables, the tight constraints of problem~\eqref{eq:integer_problem_continuous} uniquely define an optimal solution.

\paragraph{Solution strategy.}
For this class of problems the strategy corresponds to a tuple containing the index of tight constraints of the continuous reformulation~\eqref{eq:integer_problem_continuous} and the optimal value of the integer variables, \ie, ${s(\theta) = (\tightconstraints(\theta), x_{\intvars}^\star(\theta))}$.
Compared to continuous problems, we must also include the value of the integer variables to recover the optimal solution~$x^\star(\theta)$.

Given the optimal strategy, problem~\eqref{eq:integer_problem} corresponds to solving the continuous problem
\begin{equation}\label{eq:integer_reduced_problem}
\begin{array}{ll}
\text{minimize} & f(\theta, x)\\
\text{subject to} & g_i(\theta, x) \le 0, \quad \forall i \in \tightconstraints(\theta)\\
&  x_\intvars = x_{\intvars}^\star(\theta).
\end{array}
\end{equation}
Solving this problem is much less computationally demanding than~\eqref{eq:integer_problem} because it is continuous, convex and has smaller number of constraints, especially when~$n \ll m$.
As for the continuous case~\eqref{eq:continuous_reduced_problem}, the components $g_i$ that are linear in $x$ can be enforced as equalities further simplifying~\eqref{eq:integer_reduced_problem} while preserving its convexity.

Similarly to Section~\ref{sub:continuous}, in case of~\gls{MILO} and~\gls{MIQO} when the cost $f$ is linear or quadratic and the constraints $g$ are all linear, the solution of \eqref{eq:continuous_reduced_problem} corresponds to solving a linear system of equations defined by the KKT conditions~\citep[Section 10.2]{boyd2004}. This means that we can solve these problems online without needing any optimization solver~\citep{bertsimas2019}.

\subsubsection*{The knapsack problem}%
\label{ssub:knapsack_example}

Consider the knapsack problem
\begin{equation}
\begin{array}{ll}
\text{maximize} & c^\tpose x\\
\text{subject to} & a^\tpose x \le b\\
& 0 \le x \le u\\
&  x \in \integers^n,
\end{array}
\end{equation}
with $n=10$. The  decision variables are $x = (x_1, \dots, x_{10})$ indicating the quantity to pick for each item $i=1,\ldots, 10$.
We chose the cost vector $c = (0.42, 0.72, 0, 0.3, 0.15,\allowbreak 0.09, 0.19, 0.35, 0.4, 0.54)$.
The knapsack capacity is $b=5$.
The weights $a = (a_1, \dots, \allowbreak a_{10})$ and the maximum order quantity $u = (u_1, \dots, u_{10})$ are our parameters.
We assume that $a$ is  in a ball $ B(a_0, r_0)$  and $u\in  B(u_0, r_0)$  where ${u_0=a_0=(2,2,\ldots,2)}$ and $r_0=1$.

\begin{figure}
    \centering
    \begin{tikzpicture}[sibling distance=8em]
    \node[split] {$u_{2} < 2$}
        child {node[split] {$a_1 < 1.4$}
                    child {node[class3]{3}}
                    child {node[class2]{2}}
               edge from parent node[left, yshift=1ex] {True}}
        child { node[split] {$a_2 < 1.62$}
                    child {node[class4]{4}}
                    child {node[split] {$a_2 < 2.48$}
                        child {node[class1]{1}}
                        child {node[class2]{2}}
                    }
                edge from parent node[right, yshift=1ex] {False}}
    ;
\end{tikzpicture}
    \caption{Example knapsack decision strategies. \emph{Strategy 1}: $x_i = 0$ for $i\neq 2$ and $x_2 = 2$. \emph{Strategy 2}: $x_i = 0$ for $i\neq 2, 10$ and $x_2 = x_{10} = 1$.
    \emph{Strategy 3}: $x_i = 0$ for $i\neq 1, 2$ and $x_1 = 2$ and $x_2 = 1$.
    \emph{Strategy 4}: $x_i = 0$ for $i\neq 2, 10$ and $x_2 = 2$ and $x_{10} = 1$.
    }
    \label{fig:example_knapsack}
\end{figure}
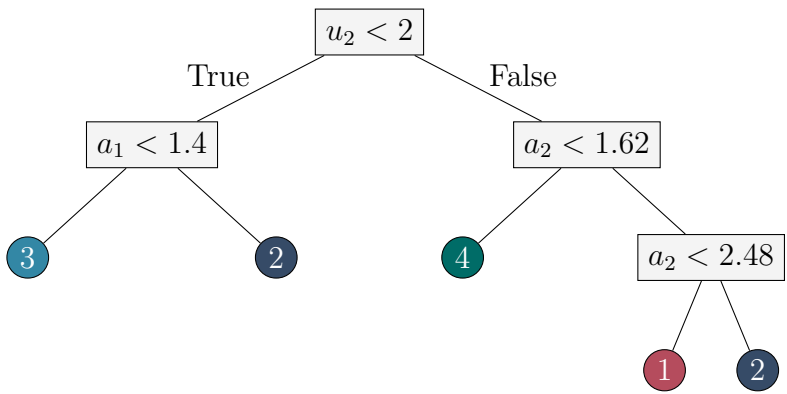

With this setup we obtain the solution strategy outlined in Figure~\ref{fig:example_knapsack}.
For this problem each strategy uniquely identifies the integer variables and is straightforward to analyze.
For instance, the left part of the tree outlines what happens if the upper bound $u_2$  is strictly less than~$2$.
Moreover, if the weight $a_1 < 1.4$, then it is easier to include $x_1$ in the knapsack solution and for this reason Strategy 3 has $x_1=2$ and $x_2=1$.
On the contrary, if $a_1 \ge 1.4$, then Strategy $2$ applies $x_1 = 0$ and $x_2=x_{10} = 1$.
Even though all these rules are simple, they capture the complexity of a hard combinatorial problem for the parameters that affect it.
Note that only a few parameters affect the strategy selection, \ie, $u_2, a_1$ and $a_2$, while the others are  not relevant for this class of problem instances
that have $a\in B(a_0, r_0)$  and $u\in  B(u_0, r_0)$.

\subsubsection*{Supplier selection}%
\label{ssub:supplier_selection_example}
Consider the problem of selecting $n=5$ suppliers in order to satisfy a known demand $d$.
Our decision variables are $x \in \{0, 1\}^n$ determining which suppliers we choose and the amount of shipments $u_i$.
For each supplier $i$ we have a per-unit cost ${c = (0.42, 0.72, 0, 0.3, 0.15)}$
and a maximum quantity to order $m = (1.09, 1.19, 1.35, 1.4 , 1.54)$, while  $\gamma = 0.1$.
We are interested in understanding the optimal strategy as a function of the
 demand $d$ and  the supplier delivery times $\tau_i$, which are the  problem parameters.
The optimization problem is as follows:
\begin{equation}
\begin{array}{ll}
\text{minimize} & \displaystyle  \sum_{i=1}^{n} c_i u_i + \gamma \max_i(\{\tau_i x_i\})\\
\text{subject to} & \displaystyle   \sum_{i=1}^{n}  u_i \ge d\\
& 0 \le u_i \le x_i m_i\\
&  x \in \{0, 1\}^n, u_i \in \reals.
\end{array}
\end{equation}

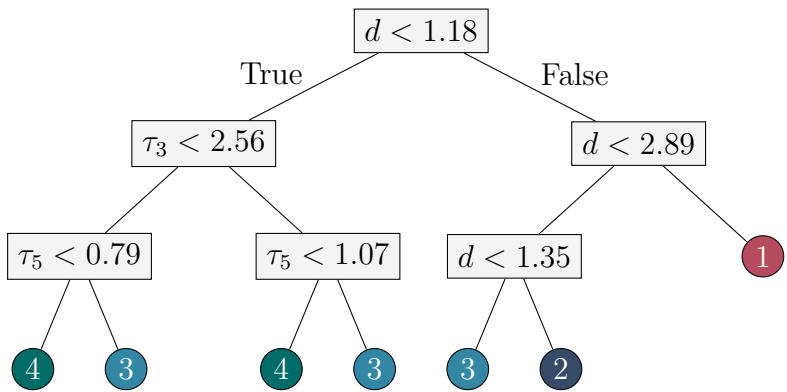
\begin{figure}
    \centering
    \begin{tikzpicture}[sibling distance=8em]
    \node[split] {$d < 1.18$}
        child {node[split] {$\tau_3 < 2.56$}
                    child {node[split]{$\tau_5 < 0.79$}
                        child {node[class4]{4}}
                        child {node[class3]{3}}
                    }
                    child {node[split]{$\tau_5 < 1.07$}
                        child {node[class4]{4}}
                        child {node[class3]{3}}
                    }
               edge from parent node[left, yshift=1ex] {True}}
        child { node[split] {$d <2.89$}
                    child {node[split]{$d < 1.35$}
                        child {node[class3]{3}}
                        child {node[class2]{2}}
                    }
                    child{node[class1]{1}}
                edge from parent node[right, yshift=1ex] {False}}
    ;
\end{tikzpicture}
    \caption{Example vendor decision strategies using \gls{OCT}.
             \emph{Strategy 1}: $x = (1, 0, 1, 0, 1)$, $u_2 = u_4 = 0$. $u_3 = m_3, u_5=m_5$ and $u_1$ is left to match demand $d$.
             \emph{Strategy 2}: $x = (0, 0, 1, 0, 1)$. $u_1 = u_2 = u_4 = 0$ and $u_3=m_3$. $u_5$ is left to match demand $d$.
             \emph{Strategy 3}: $x = (0, 0, 1, 0, 0)$. $u_1=u_2=u_4=u_5= 0$ and $u_3$ is left to match demand $d$.             \emph{Strategy 4}: $x = (0, 0, 0, 0, 1)$. $u_1=u_2=u_3=u_4=0$. $u_5$ is left to match demand $d$.
             }
    \label{fig:example_vendor}
\end{figure}

Specifically, we are interested in $d\in [1,3]$ and   $\tau = (\tau_1, \dots, \tau_5)\in B(\tau_0,r_0)$ with the center of the ball being  $\tau_0=(2, 3, 2.5, 5, 1)$
and radius $r_0=0.5$.
With these parameters we obtain the solution described in Figure~\ref{fig:example_vendor}.
Also in this case, the decision rules are simple and the solution strategy can be interpreted from the problem structure.
The strategy outputs directly the optimal choice of suppliers $x_i$ and which inequalities are tight corresponding to where we order the maximum quantity, \ie, when $u_i = m_i$.
Note that the demand constraint is always tight by construction, \ie, we want to spend the minimum to match $d$.
Therefore, we can fix $x_i$ and $u_i$ to the values given by the strategy and we can set the remaining $u_i$s to the minimum value such that $\sum_{i=1}^n u_i \ge d$.

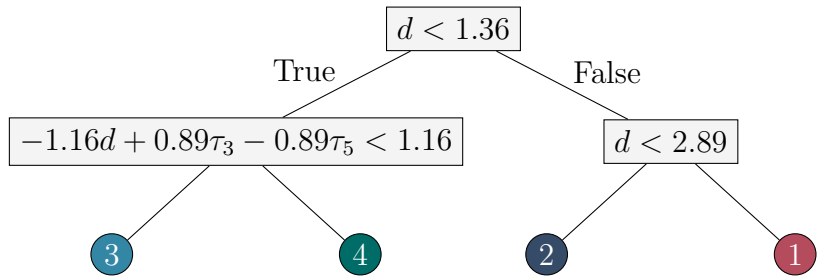
\begin{figure}
    \centering
    \begin{tikzpicture}[sibling distance=8em]
    \node[split] {$d < 1.36$}
        child { node[split] {$-1.16 d + 0.89 \tau_3 - 0.89 \tau_5 < 1.16$}
                        child {node[class3]{3}}
                        child {node[class4]{4}}
                edge from parent node[left, yshift=1ex] {True}}
        child {node[split] {$d < 2.89$}
                        child {node[class2]{2}}
                        child {node[class1]{1}}
               edge from parent node[right, yshift=1ex] {False}}
    ;
\end{tikzpicture}
    \caption{Example vendor decision strategies using \gls{OCT-H}. The strategies are the same as in Figure~\ref{fig:example_vendor} but accuracy is the same or higher and the tree depth is reduced. Smaller tree depth can sometimes help in interpreting the classification.}
    \label{fig:example_vendor_octh}
\end{figure}

In Figure~\ref{fig:example_vendor_octh},  we  show the tree that   \gls{OCT-H} gave  for this example. Even though the accuracy of the \gls{OCT} in Figure~\ref{fig:example_vendor} is 100\%, the \gls{OCT-H} depth is smaller at the cost of being less interpretable. Note that in some cases having a lower depth can make some classification tasks more interpretable despite the multiple coefficients on the hyperplanes.


\section{Machine Learning}
\label{sec:machine_learning}


In this section, we introduce  the core part of our approach: learning the  mapping from parameters to optimal  strategies.
After the training, the mapping will replace the core part of standard optimization algorithms -- the optimal strategy search -- with a multiclass classification problem where each strategy corresponds to a class label.

\subsection{Multiclass Classifier}%
\label{sub:multiclass_classifier}

We would like to solve a multiclass classification problem with i.i.d. data $(\theta_i, s_i),\; i=1,\dots,N$ where $\theta_i \in \reals^p$ are the parameters and $s_i \in \strategy$ the corresponding labels identifying the optimal strategies.
$\strategy$ represents the set of strategies of cardinality $|\strategy| = M$.

In this work we apply two supervised learning techniques for multiclass classification to compare their predictive performance: optimal classification trees (\glspl{OCT},  \glspl{OCT-H})~\citep{bertsimas2017a,bd_book} and \glsentryfullpl{NN}~\citep{bengio2009,lecun2015}.
As we mentioned earlier  \glspl{OCT}, \glspl{OCT-H} are interpretable and can be described using simple rules as in the examples in the previous section, while \glspl{NN} are not interpretable since they represent a composition  of multiple nonlinear functions.
A more detailed description of~\glspl{OCT} and \glspl{OCT-H} can be found in Appendix~\ref{sec:trees} and of \glspl{NN} can be found in Appendix~\ref{sec:neural_networks}.

\subsection{Strategies Exploration}
\label{sub:strategies_exploration}
In this section, we estimate how likely it is to find a new parameter $\theta$ whose optimal strategy does not lie among  the  ones we have already seen.
If we have already encountered most of the strategies for our problem, then it is unlikely to find unseen strategies. Therefore, we can be sure that our classification problem includes all the possible strategies (classes) arising in practice. Otherwise, we must collect more data to have a more representative classification problem.

\paragraph{Estimating the probability of finding unseen strategies.}

Given $N$ independent samples $\Theta_N = \{\theta_ 1,\ldots, \theta_N \} $ drawn from an unknown discrete distribution~$\distrib$ with the corresponding strategies
 $(s(\theta_1),\ldots, s(\theta_N))$, we find $M$ unique strategies $\strategy(\Theta_N) = \{s_1, \dots, s_M\}$.
We are interested in bounding the probability of finding unseen strategies
\begin{equation*}
\prob(s(\theta_{N+1}) \notin \strategy(\Theta_N)),
\end{equation*}
with confidence at least $1 - \beta$ where $\beta > 0$.

\paragraph{Historical background.}
This problem started from the seminal work by Turing and Good~\citep{good1953} in the context of decrypting the Enigma codes during World War II.
The Enigma machine was a German navy encryption device used for secret military communications.
Part of Enigma's encryption key was a three letter sequence (a word) selected from a book containing all the possible ones in random order.
The number of possible words was so large that it was impossible to test all the combinations with the computing power available at that time.
In order to decrypt the Enigma machine without testing all the possible words, Turing wanted to check only a subset of them while estimating that the likelihood of finding a new unseen word was low.
This is how the Good-Turing estimator was developed.
It was a fundamental step towards the Enigma machine decryption which is believed to have shortened  World War II of at least two years~\citep{turing}.
In addition, this class of estimators have become standard in a wide range of natural language processing applications.



\paragraph{Good-Turing estimator.}
Let   $N_r$ be the number of strategies that appeared exactly  $r$ times in  $(s(\theta_1),\ldots, s(\theta_N))$.
 The Good-Turing estimator for the probability of having  an  unseen strategy is given by~\citet{good1953}
\begin{equation}\label{eq:good_turing}
G = N_1 / N,
\end{equation}
which corresponds to the ratio between the number of distinct strategies that have appeared exactly once, over the total number of samples.
Despite the elegant result, Good and Turing did not provide a convergence analysis of this estimator for finite samples. Only few decades later the first theoretical work on that topic appeared in~\citep{mcallester2000}.
Using \citet{mcdiarmid1989}'s inequality the authors derived a high probability confidence interval.
That result can be directly applied to our problem with the following theorem.

\begin{theorem}[Missing strategies bound]\label{thm:missing_strategies}
The probability of encountering a parameter $\theta_{N+1}$ corresponding to an unseen strategy $s(\theta_{N+1})$ satisfies with confidence at least $1 - \beta$
\begin{equation*}
\prob(s(\theta_{N+1}) \notin \strategy(\Theta_N)) \le G + c\sqrt{(1/N)\ln(3/\beta)},
\end{equation*}
where $G$ corresponds to the Good-Turing estimator~\eqref{eq:good_turing} and $c=(2\sqrt{2} + \sqrt{3})$.
\end{theorem}

\ifpreprint
\begin{proof}
\else
\proof{Proof of Theorem~\ref{thm:missing_strategies}}
\fi
The result follows directly from~\citep[Theorem 9]{mcallester2000}.
\ifpreprint
\end{proof}
\else
\endproof
\fi

\paragraph{Exploration algorithm.}
Given the bound in Theorem~\ref{thm:missing_strategies} we  construct  Algorithm~\ref{alg:strategies_exploration} to compute the different strategies appearing in our problem.
\begin{algorithm}
  \caption{Strategies exploration}
  \label{alg:strategies_exploration}
  \begin{algorithmic}[1]
  \State {\bf given} $\epsilon, \beta, \Theta = \emptyset, \strategy = \emptyset, u = \infty$
  \For{$k = 1,\dots,$}
  \State Sample $\theta_k$ and compute $s(\theta_k)$ \Comment{Sample parameter and strategy.}
  \State $\Theta \gets \Theta \cup \{\theta_k\}$  \Comment{Update set of samples.}
  \If{$s(\theta_k) \notin \strategy$}
  \State $\strategy \gets \strategy \cup \{s(\theta_k)\}$ \Comment{Update strategy set if new strategy found}
  \EndIf
  \If{$G + c \sqrt{(1/k)\ln(3/\beta)} \le \epsilon$} \Comment{Break if bound less than $\epsilon$}
  \State {\bf break}
  \EndIf
  \EndFor
  \State \Return $k, \Theta, \strategy$
  \end{algorithmic}
\end{algorithm}
The algorithm keeps on sampling until we encounter a large enough set of distinct strategies.
It terminates when the probability of encountering a parameter with an unseen optimal strategy is less than $\epsilon$ with confidence at least $1 - \beta$.
Note that in practice, instead of sampling one point per iteration $k$, it is more convenient to sample more points to avoid having to recompute the class frequencies which becomes computationally intensive with several thousands of samples and hundreds of classes.
Algorithm~\ref{alg:strategies_exploration}   imposes no assumptions on the optimization problem nor the data distribution apart from i.i.d. samples.



\section{Machine Learning Optimizer}%
\label{sec:machine_learning_optimizer}

We implemented our method in the Python software tool MLOPT (Machine Learning Optimizer).
It is integrated with CVXPY~\citep{diamond2016} to formulate and solve the problems.
CVXPY allows the user to easily define mixed-integer convex optimization problems performing all the required reformulations needed by the optimizers while keeping track of the original constraints and variables.
This makes it ideal for identifying which constraints are tight or not at the optimal solution. The MLOPT implementation is freely available online at
\begin{center}
    \href{https://github.com/bstellato/mlopt}{\texttt{github.com/bstellato/mlopt}}.
\end{center}

We used the Gurobi Optimizer~\citep{gurobioptimization2020} to find the tight constraints because it provides a good tradeoff between solution accuracy and computation time.
Note that from Section~\ref{sec:voice}, in case of~\gls{LO},~\gls{MILO},~\gls{QO} and~\gls{MIQO} when the cost $f$ is linear or quadratic and the constraints $g$ are all linear, the  online solution corresponds to solving a linear system of equations defined by the KKT conditions~\citep[Section 10.2]{boyd2004} on the reduced subproblem. This means that we can solve those parametric optimization problems without the need to apply  any optimization  solver.

MLOPT performs the iterative sampling procedure outlined in Section~\ref{sub:strategies_exploration} to obtain the strategies required for the classification task.
We sample $5000$ new points at each iteration and compute their strategies until the Good Turing estimate is below $\eps_{\rm GT} = 0.005$.
The strategy computation is fully parallelized across samples.

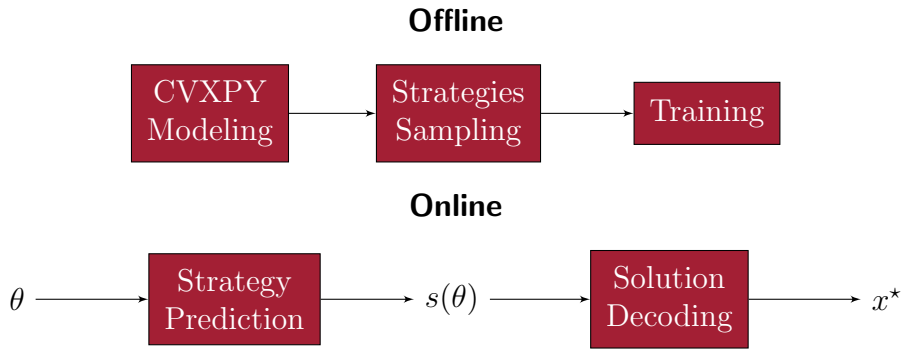
\begin{figure}
  \centering
  {\bfseries \sffamily Offline}\\[1em]
\begin{tikzpicture}[auto,>=latex', node distance=8em]
    \node [block] (cvxpy) {CVXPY\\Modeling};
    \node [block] (sampling) [right of=cvxpy] {Strategies\\Sampling};
    \path[->, draw] (cvxpy) to (sampling);
    \node [block] (training) [right of=sampling] {Training};
    \path[->, draw] (sampling) to (training);
\end{tikzpicture}\\[.5em]
  {\bfseries \sffamily Online}\\[1em]
\begin{tikzpicture}[auto,>=latex', node distance=7em]
    \node [block] (prediction) {Strategy\\Prediction};
    \node (theta)[left of=prediction] {$\theta$};
    \path[->, draw] (theta) to (prediction);
    \node (strategy) [right of=cvxpy] {$s(\theta)$};
    \path[->, draw] (prediction) to (strategy);
    \node [block] (decoding) [right of=strategy] {Solution\\Decoding};
    \path[->, draw] (strategy) to (decoding);
    \node (x)[right of=decoding] {$x^\star$};
    \path[->, draw] (decoding) to (x);
\end{tikzpicture}\\[.5em]
\caption{Algorithm implementation.}
\label{fig:block_diagram}
\end{figure}

We interfaced MLOPT to the machine learning libraries OptimalTrees.jl~\citep{bd_book} on multi-core CPUs and PyTorch~\citep{paszke2017} for both CPUs and GPUs.
In the training phase we automatically tune the predictor parameters by splitting the data points in $90\%/10\%$ training/validation.
We tune the \glspl{OCT} and \glspl{OCT-H} with maximal depth for values $5, 10, 15$ and minimum bucket size for values $1, 5, 10$.
For the \glspl{NN} we validate the stochastic gradient descent learning rate for values $0.001, 0.01, 0.1$, batch size for values $32, 128$ and number of epochs for values $50, 100$.
We described the complete algorithm in Figure~\ref{fig:block_diagram}.

\vspace{2em}

\section{Computational Benchmarks}
\label{sec:benchmarks}

In this section, we test our approach  on multiple  examples from continuous and mixed-integer optimization.
We benchmarked the predictive performance of \glspl{OCT}, \glspl{OCT-H} (see Appendix \ref{sec:trees}) and \glspl{NN} (see Appendix \ref{sec:neural_networks}) depending on the problem type and size.
We executed the numerical tests on a Dell R730 cluster with 28 Intel E5-2680 CPUs with a total of 256GB RAM and a Nvidia Tesla K80 GPU.
To facilitate the data collection, we generated the training samples from distributions on intervals or on hyperballs around specified points.
In this way the distributions are unimodal and therefore closer to the ones encountered in practice.
Note that our theoretical results do not assume anything on the distribution generating data apart from i.i.d. samples.
Therefore, we could have chosen other distributions such as multimodal ones.
In general, how to generate training data highly depends on the practical application and in many cases we can simply fit a representative distribution, either unimodal or multimodal, to the historical data we have.
We evaluated the performance metrics on 100 unseen samples drawn from the same distribution of~$\theta$.

\paragraph{Infeasibility and suboptimality.}
After the learning phase, we compared the predicted solution $\hat{x}^\star_i$ to the optimal one $x^\star_i$ obtained by solving the instances from scratch.
Given a parameter $\theta_i$, we say that the predicted solution is infeasible if the constraints are violated more than a predefined tolerance $\eps_{\rm inf} = 10^{-3}$ according to the infeasibility metric
\begin{equation*}
  p(\hat{x}^\star_i) = \|(g(\theta_i, \hat{x}^\star_i))_{+}\|_2 / r(\theta_i, \hat{x}^{\star}_i),
\end{equation*}
where $r(\theta, x)$ normalizes the violation depending on the size of the summands of $g$.
For example, if $g(\theta, x) = A(\theta)x - b$, then ${r(\theta, x) = \max (\|A(\theta)x\|_2, \|b\|_2)}$.
If the predicted solution $\hat{x}_i$ is feasible, we define its suboptimality as
\begin{equation*}
d(\hat{x}^\star_i) = (f(\hat{x}^\star_i) - f(x^\star_i))/|f(x^\star_i)|.
\end{equation*}
Note that $d(\hat{x}^\star_i) \ge 0$ by construction.
We consider a predicted solution to be accurate if it is feasible and if the suboptimality is less than the tolerance $\eps_{\rm sub} = 10^{-3}$.
For each instance we report the maximum infeasibility $ \bar{p} = \max_i p(\hat{x}^\star_i)$ and the maximum suboptimality $\bar{d} = \max_i d(\hat{x}^\star_i)$ over the test dataset.
Note that for notation ease and to simplify the ${\rm max}$ operation, if a point is infeasible we consider its suboptimality~$0$ and ignore it.

\paragraph{Predicting the optimal strategy.}
Once the learning phase is completed, the predictor outputs the  three ($3$) most likely optimal strategies
and picks the best one according to infeasibility and suboptimality after solving the associated reduced problems.

\paragraph{Runtimes.}
The training times of these examples range from a 2 to 12 hours including the time to solve the problems in parallel over the training phase and the time to train the predictors.
Even though this time can be significant, it is comparable to the usual time dedicated to train predictors for modern machine learning tasks.
We report the time ratio $t_{\rm ratio}$ between the time needed to solve the problem with Gurobi~\citep{gurobioptimization2020} compared to our method.

\paragraph{Tables notation.}
We now outline the table headings as reported in every example.
Dimensions $n$ and/or $m$ denote problem-specific sizes that vary over the benchmarks.
The total number of variables and constraints are $n_{\rm var}$ and $n_{\rm con}$ respectively.
Column ``learner'' indicates the learner type: \gls{OCT}, \gls{OCT-H} or \gls{NN}.
We indicate the Good-Turing estimator under column ``GT'' and the accuracy over the test set as ``acc [\%]''.
The column $t_{\rm ratio}$ denotes the speedup of our method over Gurobi~\citep{gurobioptimization2020}.
The maximum infeasibility appears under column $\bar{p}$ and the maximum suboptimality under column $\bar{d}$.


\subsection{Transportation Optimization}%
\label{ssub:transportation_problem}

Consider a standard transportation problem with $n$ warehouses and $m$ retail stores.
Let $x_{ij}$ be the quantity we ship from  warehouse $i$ to  store $j$.
$s_i$ denotes the supply for   warehouse  $i$ and $d_j$ the demand from  store $j$.
We define the cost of transporting a unit from warehouse $i$ to store $j$ as $c_{ij}$.
The optimization problem can be formulated as follows:
\begin{equation*}
\begin{array}{lll}
	\mbox{minimize}   & \displaystyle  \sum_{i=1} ^n \sum_{j=1}^m c_{ij} x_{ij} &\\
	\mbox{subject to} & \displaystyle \sum_{j=1}^m x_{ij} \le s_i,\quad &\forall i \\
    &  \displaystyle  \sum_{i=1} ^n  x_{ij} \ge d_j &\forall j\\
    & x_{ij} \ge 0 \quad &\forall i,j.
\end{array}
\end{equation*}
The first constraint ensures we respect the supply for each warehouse~$i$.
The second constraint enforces the sum of all the shipments to match at least the demand for store~$j$.
Finally, the quantity of shipped product has to be nonnegative.
Our learning parameters are the demands $d_j$.

\paragraph{Problem instances}

We generate the transportation problem instances by varying the number of warehouses $n$ and stores $m$.
The cost vectors $c_i$ are distributed as $\mathcal{U}(0, 5)$ and the supplies $s_i$ as $\mathcal{U}(3, 13)$.
The parameter vector $d = (d_1, \dots, d_m)$ was sampled from a uniform distribution within  the ball $B(\overline{d}, 0.75)$ with center $\overline{d} \sim \mathcal{N}(3, 1)$.

\paragraph{Results.}
The results appear in Table~\ref{tab:transportation}.
Independently from the problem instances the prediction accuracy is very high for both optimal trees and neural networks.
Note that even though we solve problems with thousands of variables and constraints, the number of strategies is always within few tens or less.
The solution time provides usually slight speedups up to almost 5 folds. Note that in the very small instances Gurobi can be faster than our method due to the delays in data exchange between predictor and our Python code.

\begin{table}
  \centering
    \caption{Transportation benchmarks.}
    \label{tab:transportation}
    \begin{adjustbox}{max width=\textwidth}
        \begin{tabular}{
                l
                S[table-format=2, round-precision=0]
                S[table-format=4, round-precision=0]
                S[table-format=4, round-precision=0]
                l
                S[table-format=5, round-precision=0]
                S[table-format = 1.2e1,
                  scientific-notation=true,
                  retain-zero-exponent=true,
                  round-integer-to-decimal]
                S[table-format=2, round-precision=0]
                S[table-format=3.2, round-integer-to-decimal]
                S[table-format=3.2, round-integer-to-decimal]
                S[table-format = 1.2e2,
                  scientific-notation=true,
                  retain-zero-exponent=true,
                  round-integer-to-decimal]
                S[table-format = 1.2e2,
                  scientific-notation=true,
                  retain-zero-exponent=true,
                  round-integer-to-decimal]
                }
\toprule
$n$ & $m$ & $n_{\rm var}$ & $n_{\rm con}$ & learner & $N$ & GT & $|\strategy|$ & ${{\rm acc}\;[\%]}$ & $t_{\rm ratio}$ & ${\bar{p}}$ & ${\bar{d}}$ \\
\midrule
\begin{filecontents*}{./data/benchmarks/transportation/transportation_general.csv}
,accuracy,avg_infeas,avg_subopt,avg_time_improv,good_turing,good_turing_smooth,learner,m,max_infeas,max_subopt,max_time_improv,mean_pred_time_pred,mean_solve_time_pred,mean_time_full,mean_time_pred,n,n_best,n_constr,n_correct,n_infeas,n_strategies,n_test,n_theta,n_train,n_var,predictor,problem,std_infeas,std_pred_time_pred,std_solve_time_pred,std_subopt,std_time_full,std_time_pred
0,100.0,6.342080052344551e-07,-1.7425487476558042e-07,0.8119282796569325,0.0001,0.00020296000000000018,pytorch,20.0,6.342080051374821e-05,4.403767040416214e-16,1.0864057192944012,1.3656616210937493e-05,0.0010054349899291993,0.0008274292945861816,0.001019091606140137,20.0,10.0,440.0,100.0,0.0,11.0,100.0,20.0,10000.0,400.0,NN,transportation,6.310289976232709e-06,6.776263578034403e-21,0.0002221932133727086,1.7338141123253697e-06,0.00017253518395713384,0.00022219321337270853
1,100.0,6.342080052344551e-07,-1.7425487476558042e-07,0.04272242483561129,0.0001,0.00020296000000000018,optimaltree,20.0,6.342080051374821e-05,4.403767040416214e-16,0.08894738903306966,0.023570799827575678,0.0009833574295043945,0.0010490131378173828,0.024554157257080073,20.0,10.0,440.0,100.0,0.0,11.0,100.0,20.0,10000.0,400.0,OCT,transportation,6.310289976232709e-06,6.938893903907228e-18,0.0002314320182370296,1.7338141123253697e-06,0.000258948273200344,0.0002314320182370296
2,100.0,6.342080052344551e-07,-1.7425487476558042e-07,1.2162547585235748,0.0001,0.00020296000000000018,optimaltree,20.0,6.342080051374821e-05,4.403767040416214e-16,1.5927858439201452,4.9591064453125e-05,0.0007495594024658203,0.0009719705581665039,0.0007991504669189453,20.0,10.0,440.0,100.0,0.0,11.0,100.0,20.0,10000.0,400.0,OCT-H,transportation,6.310289976232709e-06,0.0,4.5244199189066475e-05,1.7338141123253697e-06,0.0002266989782403518,4.5244199189066475e-05
0,100.0,7.240135579737208e-17,-7.256169511536382e-18,0.743075892192386,0.0001,0.00020296000000000018,pytorch,10.0,1.4184731607937946e-16,2.761345960927436e-16,0.816097867718516,8.64267349243164e-06,0.0007599043846130371,0.0005710887908935547,0.0007685470581054687,20.0,6.0,230.0,100.0,0.0,6.0,100.0,10.0,10000.0,200.0,NN,transportation,6.557711303981488e-17,0.0,6.089445614999618e-05,5.381901055303088e-17,7.436063264446375e-05,6.089445614999618e-05
1,100.0,7.240135579737208e-17,-7.256169511536382e-18,0.8757981874413376,0.0001,0.00020296000000000018,optimaltree,10.0,1.4184731607937946e-16,2.761345960927436e-16,1.4190553469326137,4.9951076507568344e-05,0.000724802017211914,0.0006785273551940918,0.0007747530937194827,20.0,6.0,230.0,100.0,0.0,6.0,100.0,10.0,10000.0,200.0,OCT,transportation,6.557711303981488e-17,1.3552527156068805e-20,3.611425183537422e-05,5.381901055303088e-17,0.00015647595118019355,3.6114251835374216e-05
2,100.0,7.240135579737208e-17,-7.256169511536382e-18,1.3575516408834845,0.0001,0.00020296000000000018,optimaltree,10.0,1.4184731607937946e-16,2.761345960927436e-16,2.155491099357334,6.237506866455076e-05,0.0004175472259521484,0.0006515192985534668,0.000479922294616699,20.0,6.0,230.0,100.0,0.0,6.0,100.0,10.0,10000.0,200.0,OCT-H,transportation,6.557711303981488e-17,2.710505431213761e-20,2.6833010251951813e-05,5.381901055303088e-17,0.0001823343554682562,2.683301025195181e-05
0,100.0,1.137254757355475e-16,-3.802561184478632e-18,1.9537342131713495,0.0002,0.00030592000000000016,pytorch,40.0,3.0036332491579103e-16,1.9811832594476118e-16,2.696297278665022,1.569032669067383e-05,0.0016210293769836426,0.0031977152824401857,0.0016367197036743168,40.0,9.0,1680.0,100.0,0.0,9.0,100.0,40.0,10000.0,1600.0,NN,transportation,6.822245287427427e-17,0.0,0.0002227666213052391,9.123974003036528e-17,0.0004775648354642239,0.0002227666213052391
1,100.0,1.137254757355475e-16,-3.802561184478632e-18,1.7356990546638376,0.0002,0.00030592000000000016,optimaltree,40.0,3.0036332491579103e-16,1.9811832594476118e-16,2.0162454873646207,8.869171142578125e-05,0.0017304635047912598,0.0031575059890747072,0.001819155216217041,40.0,9.0,1680.0,100.0,0.0,9.0,100.0,40.0,10000.0,1600.0,OCT,transportation,6.822245287427427e-17,0.0,9.785672732677236e-05,9.123974003036528e-17,0.000434418456189302,9.785672732677236e-05
2,100.0,1.137254757355475e-16,-3.802561184478632e-18,1.8140458336345535,0.0002,0.00030592000000000016,optimaltree,40.0,3.0036332491579103e-16,1.9811832594476118e-16,2.2563363302605346,8.934020996093753e-05,0.0017245364189147948,0.0032904553413391113,0.0018138766288757322,40.0,9.0,1680.0,100.0,0.0,9.0,100.0,40.0,10000.0,1600.0,OCT-H,transportation,6.822245287427427e-17,2.710505431213761e-20,7.268791571166037e-05,9.123974003036528e-17,0.0005655268674470461,7.268791571166037e-05
0,100.0,1.131628409945679e-16,4.066420790359955e-18,4.978899249212885,0.0,0.00010000000000000018,pytorch,20.0,1.408265175379916e-16,2.1081295767982203e-16,4.824884751694306,8.02516937255859e-06,0.0003069925308227539,0.001568441390991211,0.0003150177001953125,40.0,3.0,860.0,100.0,0.0,3.0,100.0,20.0,10000.0,800.0,NN,transportation,4.157888195577416e-17,3.3881317890172014e-21,4.238290814757885e-05,4.1022839412475416e-17,0.00018988547886146596,4.238290814757885e-05
1,100.0,1.131628409945679e-16,4.066420790359955e-18,3.5270116410453998,0.0,0.00010000000000000018,optimaltree,20.0,1.408265175379916e-16,2.1081295767982203e-16,3.807706198290972,6.034135818481446e-05,0.0004903888702392578,0.0019424319267272948,0.0005507302284240722,40.0,3.0,860.0,100.0,0.0,3.0,100.0,20.0,10000.0,800.0,OCT,transportation,4.157888195577416e-17,6.776263578034403e-21,3.064986925586074e-05,4.1022839412475416e-17,0.0002259932633376606,3.064986925586074e-05
2,100.0,1.131628409945679e-16,4.066420790359955e-18,3.737755774078534,0.0,0.00010000000000000018,optimaltree,20.0,1.408265175379916e-16,2.1081295767982203e-16,4.105294974681524,5.3045749664306656e-05,0.00047186851501464846,0.001962001323699951,0.0005249142646789548,40.0,3.0,860.0,100.0,0.0,3.0,100.0,20.0,10000.0,800.0,OCT-H,transportation,4.157888195577416e-17,1.3552527156068805e-20,3.13933035726056e-05,4.1022839412475416e-17,0.0002483512643633493,3.13933035726056e-05
0,99.0,0.00018650909316837955,-8.677366912489569e-09,2.044370950903137,0.0015,0.00160084,pytorch,60.0,0.018591889163771996,2.4561914493787986e-16,2.2755569121597508,1.0046958923339844e-05,0.003324141502380371,0.006816318035125732,0.0033341884613037097,60.0,10.0,3720.0,99.0,1.0,43.0,100.0,60.0,10000.0,3600.0,NN,transportation,0.001849819617889809,0.0,0.00019164487420929115,8.590154982600873e-08,0.0006359617448229711,0.0001916448742092912
1,95.0,0.0007444181024967428,-9.042729729857553e-09,1.8025150053890218,0.0015,0.00160084,optimaltree,60.0,0.027200046680806184,0.00015532382649657736,2.5211106717860514,0.0001679730415344238,0.0035416793823242186,0.006686704158782959,0.0037096524238586425,60.0,10.0,3720.0,95.0,5.0,43.0,100.0,60.0,10000.0,3600.0,OCT,transportation,0.0036486733890343775,2.710505431213761e-20,0.00026229852113770157,8.767251754794448e-08,0.0007757936441967979,0.00026229852113770146
2,100.0,5.170406355182213e-07,-6.373452094175094e-08,2.109531497697329,0.0015,0.00160084,optimaltree,60.0,5.1704063539725705e-05,2.4561914493787986e-16,2.6671881570705573,0.00016964673995971684,0.0032533860206604005,0.007220995426177978,0.003423032760620117,60.0,10.0,3720.0,100.0,0.0,43.0,100.0,60.0,10000.0,3600.0,OCT-H,transportation,5.144489366926793e-06,5.421010862427522e-20,0.0002594256644057841,6.341504764831421e-07,0.0011622897575448916,0.0002594256644057841
0,99.0,1.6040083720698816e-05,-2.1059150543884237e-17,1.6999743438324562,0.0008,0.0009137800000000001,pytorch,30.0,0.0016040083720644787,1.9700657991029576e-16,2.388313421149859,1.9202232360839842e-05,0.0021367359161376955,0.0036650395393371584,0.0021559381484985356,60.0,10.0,1890.0,99.0,1.0,30.0,100.0,30.0,10000.0,1800.0,NN,transportation,0.00015959681792179436,3.3881317890172014e-21,0.0001528165500220738,1.1409534781351955e-16,0.00038826245091936175,0.00015281655002207387
1,100.0,5.402645910528175e-17,-2.0848559038445397e-17,1.4787913130541133,0.0008,0.0009137800000000001,optimaltree,30.0,1.1903425381757013e-16,1.9700657991029576e-16,1.231452573651747,0.00013025522232055665,0.0022343492507934572,0.003496756553649902,0.002364604473114014,60.0,10.0,1890.0,100.0,0.0,30.0,100.0,30.0,10000.0,1800.0,OCT,transportation,5.387035654319484e-17,0.0,0.00024756380362571406,1.1354277360622848e-16,0.00044152059333055693,0.0002475638036257141
2,100.0,5.402645910528175e-17,-2.0848559038445397e-17,1.5545680417100103,0.0008,0.0009137800000000001,optimaltree,30.0,1.1903425381757013e-16,1.9700657991029576e-16,2.1642460799542977,0.00014233112335205079,0.0021290183067321777,0.0035309672355651857,0.0022713494300842275,60.0,10.0,1890.0,100.0,0.0,30.0,100.0,30.0,10000.0,1800.0,OCT-H,transportation,5.387035654319484e-17,0.0,0.00010900843886274822,1.1354277360622848e-16,0.0005430297474760489,0.00010900843886274822
0,100.0,1.5732632222945104e-16,-2.219757063074307e-18,2.0445132924801244,0.0011,0.0012068200000000004,pytorch,80.0,1.9884945748529271e-16,2.1995615585506774e-16,2.4672128719935595,1.7852783203125e-05,0.005424394607543946,0.011126747131347656,0.005442247390747069,80.0,10.0,6560.0,100.0,0.0,36.0,100.0,80.0,10000.0,6400.0,NN,transportation,1.3895655148365712e-17,0.0,0.0003122266325717938,9.5469663266235e-17,0.0010213854198918592,0.0003122266325717938
1,100.0,1.5732632222945104e-16,-2.219757063074307e-18,2.1150678658203055,0.0011,0.0012068200000000004,optimaltree,80.0,1.9884945748529271e-16,2.1995615585506774e-16,1.8587049277039154,0.00016175746917724616,0.005460939407348633,0.011892385482788086,0.005622696876525878,80.0,10.0,6560.0,100.0,0.0,36.0,100.0,80.0,10000.0,6400.0,OCT,transportation,1.3895655148365712e-17,5.421010862427522e-20,0.0006406436922757174,9.5469663266235e-17,0.002102028672390334,0.0006406436922757174
2,100.0,1.5732632222945104e-16,-2.219757063074307e-18,2.4262071082132497,0.0011,0.0012068200000000004,optimaltree,80.0,1.9884945748529271e-16,2.1995615585506774e-16,2.934232911850449,0.00013578653335571283,0.004245576858520508,0.0106300950050354,0.004381363391876221,80.0,10.0,6560.0,100.0,0.0,36.0,100.0,80.0,10000.0,6400.0,OCT-H,transportation,1.3895655148365712e-17,5.421010862427522e-20,0.0003488198551514174,9.5469663266235e-17,0.0022712802197956,0.0003488198551514174
0,100.0,1.3615751443209418e-16,2.0580439386885848e-17,2.3992539631972276,0.0001,0.00020296000000000018,pytorch,40.0,2.1041052813001341e-16,3.2318731095032843e-16,2.096687954405606,1.743316650390625e-05,0.002399871349334717,0.005799727439880371,0.0024173045158386226,80.0,10.0,3320.0,100.0,0.0,14.0,100.0,40.0,10000.0,3200.0,NN,transportation,6.282040613032422e-17,0.0,0.00015483615965349905,8.623356253827203e-17,0.0007546346854115568,0.00015483615965349907
1,100.0,1.3615751443209418e-16,2.0580439386885848e-17,2.1268876341265526,0.0001,0.00020296000000000018,optimaltree,40.0,2.1041052813001341e-16,3.2318731095032843e-16,2.4776980170501903,6.395101547241211e-05,0.002812814712524414,0.006118557453155518,0.0028767657279968255,80.0,10.0,3320.0,100.0,0.0,14.0,100.0,40.0,10000.0,3200.0,OCT,transportation,6.282040613032422e-17,0.0,0.00016330570684970066,8.623356253827203e-17,0.0007664620817454344,0.00016330570684970066
2,100.0,1.3615751443209418e-16,2.0580439386885848e-17,2.3839053077507923,0.0001,0.00020296000000000018,optimaltree,40.0,2.1041052813001341e-16,3.2318731095032843e-16,2.936188765076971,5.615234375000001e-05,0.0023678088188171386,0.0057784938812255856,0.0024239611625671395,80.0,10.0,3320.0,100.0,0.0,14.0,100.0,40.0,10000.0,3200.0,OCT-H,transportation,6.282040613032422e-17,6.776263578034403e-21,0.00019627544118625795,8.623356253827203e-17,0.0012365939223351193,0.00019627544118625795
\end{filecontents*}
\csvreader[head to column names, late after line=\\]{./data/benchmarks/transportation/transportation_general.csv}{
predictor=\predictor,
accuracy=\accuracy,
n_strategies=\nstrategies,
n_infeas=\ninfeas,
avg_infeas=\avginfeas,
avg_subopt=\avgsubopt,
max_infeas=\maxinfeas,
max_subopt=\maxsubopt,
avg_time_improv=\avgtimeimprov,
max_time_improv=\maxtimeimprov,
n=\n,
m=\m,
n_var=\nvar,
n_constr=\nconstr,
n_train=\ntrain,
good_turing=\goodturing,
}
{\tablenum[table-format=2, round-precision=0]{\n} & \m & \nvar & \nconstr & \predictor & \ntrain & \goodturing &\nstrategies &\accuracy & \avgtimeimprov & \maxinfeas & \maxsubopt}
     \bottomrule
\end{tabular}
\end{adjustbox}
\end{table}

\subsection{Portfolio Optimization}%
\label{ssub:portfolio}

Consider the problem of allocating assets to minimize the risk adjusted return considered in \cite{markowitz1952} who formulated it as a \gls{QO}
\begin{equation*}
\begin{array}{ll}
	\mbox{maximize}   & \mu^\tpose x - \gamma (x^\tpose \Sigma x) \\
	\mbox{subject to} & \ones^\tpose x = 1 \\
	& x \ge 0,
\end{array}
\end{equation*}
where the variable $x \in \reals^{n}$ represents the investments to make.

Our learning parameter is the vector of expected returns $\mu  \in \reals^{n}$.
In addition we denote the risk aversion coefficient as $\gamma > 0$, and the covariance matrix for the risk model as $\Sigma\in\symm_{+}^n$. $\Sigma$ is
 assumed to be
\begin{equation*}
    \Sigma = FF^\tpose + D,
\end{equation*}
where $F\in \reals^{n\times p}$ is   the factor loading matrix and $D\in \reals^{n\times n}$ is a diagonal matrix describing the asset-specific risk.

\paragraph{Problem instances.}
We generated portfolio instances for different number of factors~$p$ and assets~$n$.
We chose the elements of $F$ with $50\;\%$ nonzeros and $F_{ij} \sim \mathcal{N}(0,1)$.
The diagonal matrix $D$ has elements $D_{ii} \sim \mathcal{U}(0, \sqrt{p})$.
The return vector parameters were randomly sampled  from a uniform distribution within  the ball $B(\overline{\mu}, 0.15)$ with $\overline{\mu} \sim \mathcal{N}(0, 1)$.
The risk-aversion coefficient is $\gamma = 1$.

\paragraph{Results.}
Results are shown in Table~\ref{tab:portfolio}.
The prediction accuracy   for \glspl{OCT},  \glspl{OCT-H} and \glspl{NN} is 100\%, the number of strategies is less than 15, and the infeasibility and suboptimality  are very low. The solution times improvement is not significant due to the delays of passing data to the predictor and the fact that Gurobi is fast for these problems.

\begin{table}
  \centering
    \caption{Continuous portfolio benchmarks.}
    \label{tab:portfolio}
    \begin{adjustbox}{max width=\textwidth}
        \begin{tabular}{
                l
                S[table-format=2, round-precision=0]
                S[table-format=3, round-precision=0]
                l
                S[table-format=5, round-precision=0]
                S[table-format = 1.2e1,
                  scientific-notation=true,
                  retain-zero-exponent=true,
                  round-integer-to-decimal]
                S[table-format=2, round-precision=0]
                S[table-format=3.2, round-integer-to-decimal]
                S[table-format=1.2, round-integer-to-decimal]
                S[table-format = 1.2e2,
                  scientific-notation=true,
                  retain-zero-exponent=true,
                  round-integer-to-decimal]
                S[table-format = 1.2e2,
                  scientific-notation=true,
                  retain-zero-exponent=true,
                  round-integer-to-decimal]
                }
\toprule
$n$ & $p$ & $n_{\rm con}$ & learner & $N$ & ${\rm GT}$ & $|\strategy|$ & ${\rm acc}\;[\%]$ & $t_{\rm ratio}$  & ${\bar{p}}$ & ${\bar{d}}$ \\
\midrule
\begin{filecontents*}{./data/benchmarks/portfolio/portfolio_general.csv}
,accuracy,avg_infeas,avg_subopt,avg_time_improv,good_turing,good_turing_smooth,learner,max_infeas,max_subopt,max_time_improv,mean_pred_time_pred,mean_solve_time_pred,mean_time_full,mean_time_pred,n_best,n_constr,n_correct,n_infeas,n_strategies,n_test,n_theta,n_train,n_var,p,predictor,problem,std_infeas,std_pred_time_pred,std_solve_time_pred,std_subopt,std_time_full,std_time_pred
0,100.0,1.654232306691483e-16,2.4598472599878665e-17,1.0138560046386926,0.0004,0.0004979800000000002,pytorch,6.661338147750934e-16,7.062827727289685e-16,1.2165561526411677,1.5845298767089842e-05,0.0018098044395446773,0.0018509459495544434,0.001825649738311768,10.0,101.0,100.0,0.0,12.0,100.0,100.0,10000.0,100.0,10.0,NN,portfolio,1.6946455303554792e-16,3.3881317890172006e-21,0.0001900601991466983,2.1231270692118541e-16,0.00044231970902890264,0.00019006019914669824
1,100.0,1.654232306691483e-16,2.4598472599878665e-17,0.08885373270854133,0.0004,0.0004979800000000002,optimaltree,6.661338147750934e-16,7.062827727289685e-16,0.12748045920803874,0.023621094226837173,0.0018283557891845704,0.0022612786293029785,0.02544945001602174,10.0,101.0,100.0,0.0,12.0,100.0,100.0,10000.0,100.0,10.0,OCT,portfolio,1.6946455303554792e-16,1.3877787807814454e-17,0.0001295844474022095,2.1231270692118541e-16,0.00046863302459810647,0.0001295844474022095
2,100.0,1.654232306691483e-16,2.4598472599878665e-17,1.0749975997638876,0.0004,0.0004979800000000002,optimaltree,6.661338147750934e-16,7.062827727289685e-16,1.6013880118610468,9.393453598022456e-05,0.0019175243377685548,0.002162313461303711,0.00201145887374878,10.0,101.0,100.0,0.0,12.0,100.0,100.0,10000.0,100.0,10.0,OCT-H,portfolio,1.6946455303554792e-16,4.0657581468206416e-20,0.0002521424241391435,2.1231270692118541e-16,0.000733848425988955,0.0002521424241391435
0,100.0,9.992007221626408e-17,2.7739710583632783e-17,0.8783905495553647,0.0004,0.0004979800000000002,pytorch,4.440892098500626e-16,5.539909432625092e-16,1.0647617061215977,2.3913383483886705e-05,0.009934601783752441,0.00874746561050415,0.009958515167236323,10.0,201.0,100.0,0.0,11.0,100.0,200.0,10000.0,200.0,20.0,NN,portfolio,1.3551999095955062e-16,1.3552527156068802e-20,0.00047128886233627323,2.2346553726929006e-16,0.001688571840818083,0.00047128886233627334
1,100.0,9.769962616701378e-17,3.83682694895117e-09,0.8132648596309525,0.0004,0.0004979800000000002,optimaltree,4.440892098500626e-16,3.836826918437501e-07,0.9332980435023394,7.446527481079103e-05,0.009651641845703123,0.007909901142120361,0.00972610712051392,10.0,201.0,100.0,0.0,11.0,100.0,200.0,10000.0,200.0,20.0,OCT,portfolio,1.3531974638774172e-16,1.3552527156068802e-20,0.0005227003731056769,3.817594581891146e-08,0.0015128285781514992,0.0005227003731056769
2,100.0,9.992007221626408e-17,2.7739710583632783e-17,0.7879980137605241,0.0004,0.0004979800000000002,optimaltree,4.440892098500626e-16,5.539909432625092e-16,1.3251297697833642,4.9383640289306673e-05,0.009543826580047607,0.0075594305992126454,0.009593210220336913,10.0,201.0,100.0,0.0,11.0,100.0,200.0,10000.0,200.0,20.0,OCT-H,portfolio,1.3551999095955062e-16,2.7105054312137605e-20,0.0004690021204961757,2.2346553726929006e-16,0.001674872627154261,0.0004690021204961757
0,100.0,3.382648696880168e-07,-1.0091619948400078e-12,2.8666863108526344,0.0001,0.00020098000000000018,pytorch,3.382648695758843e-05,6.960303022906644e-16,3.124961719831415,1.5528202056884765e-05,0.007735874652862548,0.022220840454101564,0.007751402854919431,4.0,301.0,100.0,0.0,4.0,100.0,300.0,10000.0,300.0,30.0,NN,portfolio,3.365692956413872e-06,0.0,0.000661578878588018,1.004095398373716e-11,0.003080506564301544,0.0006615788785880178
1,100.0,3.382648696880168e-07,-1.0091619948400078e-12,2.5476091377600905,0.0001,0.00020098000000000018,optimaltree,3.382648695758843e-05,6.960303022906644e-16,3.0204082445563736,4.0311813354492204e-05,0.0077100348472595205,0.019744853973388668,0.0077503466606140175,4.0,301.0,100.0,0.0,4.0,100.0,300.0,10000.0,300.0,30.0,OCT,portfolio,3.365692956413872e-06,6.776263578034403e-21,0.0006658744941214018,1.004095398373716e-11,0.003700769771283992,0.0006658744941214018
2,100.0,3.382648696880168e-07,-1.0091619948400078e-12,2.461744987279796,0.0001,0.00020098000000000018,optimaltree,3.382648695758843e-05,6.960303022906644e-16,3.905469656367532,3.5452842712402355e-05,0.006568748950958252,0.016257860660552983,0.0066042017936706575,4.0,301.0,100.0,0.0,4.0,100.0,300.0,10000.0,300.0,30.0,OCT-H,portfolio,3.365692956413872e-06,1.3552527156068802e-20,0.000577971725049501,1.004095398373716e-11,0.0039478643578307535,0.0005779717250495007
0,100.0,1.475364763656361e-07,-9.345732212287949e-14,1.2879392960444358,0.0,0.00010396000000000019,pytorch,1.4753647625239334e-05,5.259838452574184e-16,1.1191323246445826,2.2745132446289057e-05,0.0352367115020752,0.0454120397567749,0.03525945663452148,10.0,401.0,100.0,0.0,14.0,100.0,400.0,10000.0,400.0,40.0,NN,portfolio,1.467969403849724e-06,6.776263578034403e-21,0.0020379482196251314,9.296373511541495e-13,0.004324587482638628,0.002037948219625131
1,100.0,1.475364763656361e-07,-9.345732212287949e-14,0.7230700209436094,0.0,0.00010396000000000019,optimaltree,1.4753647625239334e-05,5.259838452574184e-16,0.9172742872262246,0.023253557682037358,0.036866748332977296,0.043471190929412845,0.06012030601501465,10.0,401.0,100.0,0.0,14.0,100.0,400.0,10000.0,400.0,40.0,OCT,portfolio,1.467969403849724e-06,3.469446951953614e-18,0.002172473809916594,9.296373511541495e-13,0.004690433718176738,0.002172473809916594
2,100.0,1.4753647636785655e-07,1.8199994513099833e-09,1.1990388009048147,0.0,0.00010396000000000019,optimaltree,1.4753647625239334e-05,1.82009291388845e-07,2.309843706431594,8.117675781250003e-05,0.03641039371490479,0.04375480890274048,0.03649157047271729,10.0,401.0,100.0,0.0,14.0,100.0,400.0,10000.0,400.0,40.0,OCT-H,portfolio,1.467969403849501e-06,2.710505431213761e-20,0.0019695830468027965,1.8109705254058837e-08,0.009677360324459814,0.0019695830468027965
\end{filecontents*}
\csvreader[head to column names, late after line=\\]{./data/benchmarks/portfolio/portfolio_general.csv}{
predictor=\predictor,
accuracy=\accuracy,
n_strategies=\nstrategies,
n_infeas=\ninfeas,
avg_infeas=\avginfeas,
avg_subopt=\avgsubopt,
max_infeas=\maxinfeas,
max_subopt=\maxsubopt,
avg_time_improv=\avgtimeimprov,
max_time_improv=\maxtimeimprov,
p=\p,
n_var=\nvar,
n_constr=\nconstr,
n_train=\ntrain,
good_turing=\goodturing,
}
{\tablenum[table-format=3, round-precision=0]{\nvar} & \p & \nconstr &\predictor & \ntrain & \goodturing &\nstrategies &\accuracy & \avgtimeimprov & \maxinfeas & \maxsubopt } 
     \bottomrule
\end{tabular}
\end{adjustbox}
\end{table}


\subsection{Facility location}%
\label{ssub:facility_location}
Let $i \in I$ the set of possible locations for facilities such as factories or warehouses.
We denote as $j \in J$ the set of delivery locations.
The cost of transporting a unit of goods from facility $i$ to location $j$ is $c_{ij}$.
The construction cost for building facility $i$ is $f_i$.
For each location $j$, we define the demand~$d_j$.
The capacity of facility $i$ is $s_i$.
The main goal of this problem is to find the best tradeoff between transportation costs versus construction costs.
We can write the model as a \gls{MILO}
\begin{equation}\label{eq:facility_location}
\begin{array}{ll}
    \text{minimize} &  \displaystyle \sum_{i\in I} \sum_{j\in J} c_{ij} x_{ij} + \sum_{i\in I} f_i y_i\\
\text{subject to} & \displaystyle  \sum_{i\in I } x_{ij} \ge d_{j},\quad \forall j \in J\\
& \displaystyle  \sum_{j\in J } x_{ij} \le s_i y_i, \quad \forall i \in I\\
&  x_{ij} \ge 0\\
&  y_i \in \{0, 1\}.
\end{array}
\end{equation}
The decision variables $x_{ij}$ describe the amount of goods sent from facility $i$ to location $j$.
The binary variables $y_i$ determine if we build facility $i$ or not.

\paragraph{Problem instances.}
We generated the facility location instances for different values of facilities $n$ and warehouses $m$.
The costs $c_{ij}$ are sampled from a uniform distribution $\mathcal{U}(0, 1)$ and $f_i$ from $\mathcal{U}(0, 10)$.
We chose capacities $s_i$ as $\mathcal{U}(8, 18)$ to ensure the problem is always feasible.
The demands parameters $d = (d_1, \dots, d_m)$ were sampled from a uniform distributions within the  ball  $B(\overline{d}, 0.25)$ with $\overline{d} \sim \mathcal{N}(3, 1)$.

\begin{figure}
    \centering
    \begin{tikzpicture}[level 1/.style={sibling distance=18em},
         level 2/.style={sibling distance=12em},
         level 3/.style={sibling distance=8em}]
        \node[split] {$d_{8} < 3.35$}
            child { node[split] {$d_{18} < 3.04$}
                    child {node[class1]{1}}
                    child {node[split] {$d_{30} < 2.19$}
                        child{node[class2]{2}}
                        child { node[split] {$d_{18} < 3.04$}
                            child {node[class2]{2}}
                            child {node[class1]{1}}
                        }
                }
            edge from parent node[left, yshift=1ex] {False}
            }
            child { node[class1]{1}
            edge from parent node[right, yshift=1ex] {True}
            }
    ;
\end{tikzpicture}
\caption{Facility location benchmark \gls{OCT} for $n$ = 60 and $m$ = 30.}
    \label{fig:benchmark_facility}
\end{figure}
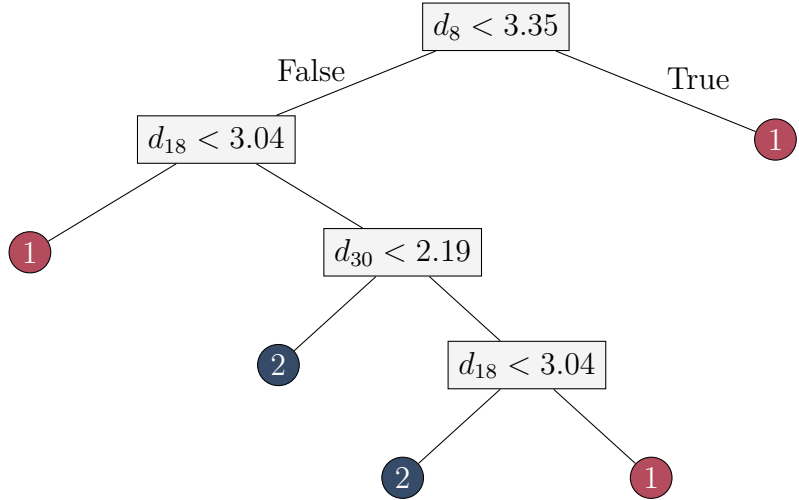

\paragraph{Results.}
Table~\ref{tab:facility} shows the benchmark examples.
All prediction methods showed very high accuracy.
There is a significant improvement in computation time compared to Gurobi.
To illustrate an example of the \gls{OCT} interpretability, we reported the resulting tree for $n=60$ and $m=30$ in Figure~\ref{fig:benchmark_facility}.
Even though for space concerns we cannot directly report the strategies in terms of tight constraints and integer variables since the problem involves $1860$ variables and $2010$ constraints, the tree is remarkably simple. In every strategy all the demand constraints are tight since we are trying to minimize the cost while satisfying the demand.
In addition, the strategy illustrates which facilities hit the maximum capacity $s_i$ which could mean that further facilities might be needed in that area. Finally, with the values of binary variables $y_i$, the strategy tells us which facilities we should open depending on the demand.

\begin{table}
  \centering
    \caption{Facility location benchmarks.}
    \label{tab:facility}
    \begin{adjustbox}{max width=\textwidth}
\begin{tabular}{l
                S[table-format=2, round-precision=0]
                S[table-format=4, round-precision=0]
                S[table-format=4, round-precision=0]
                l
                S[table-format=5, round-precision=0]
                S[table-format = 1.2e1,
                  scientific-notation=true,
                  retain-zero-exponent=true,
                  round-integer-to-decimal]
                S[table-format=1, round-precision=0]
                S[table-format=3.2, round-integer-to-decimal]
                S[table-format=3.2, round-integer-to-decimal]
                S[table-format = 1.2e2,
                  scientific-notation=true,
                  retain-zero-exponent=true,
                  round-integer-to-decimal]
                S[table-format = 1.2e2,
                  scientific-notation=true,
                  retain-zero-exponent=true,
                  round-integer-to-decimal]
                }
\toprule
$n$ & $m$ & $n_{\rm var}$ & $n_{\rm con}$ & learner & $N$ & ${\rm GT}$ & $|\strategy|$ & ${\rm acc}\;[\%]$ & $t_{\rm ratio}$ & ${\bar{p}}$ & ${\bar{d}}$ \\
\midrule
\begin{filecontents*}{./data/benchmarks/facility/facility_general.csv}
,accuracy,avg_infeas,avg_subopt,avg_time_improv,good_turing,good_turing_smooth,learner,m,max_infeas,max_subopt,max_time_improv,mean_pred_time_pred,mean_solve_time_pred,mean_time_full,mean_time_pred,n,n_best,n_constr,n_correct,n_infeas,n_strategies,n_test,n_theta,n_train,n_var,predictor,problem,std_infeas,std_pred_time_pred,std_solve_time_pred,std_subopt,std_time_full,std_time_pred
0,100.0,3.7317930012595747e-17,1.64464618860614e-18,537.6686771604634,0.0,0.00010000000000000018,pytorch,10.0,1.0638049822737197e-16,1.64464618860614e-16,1482.1299795551704,8.118152618408206e-06,0.00013774394989013672,0.07842548370361327,0.00014586210250854497,20.0,1.0,270.0,100.0,0.0,1.0,100.0,10.0,10000.0,220.0,NN,facility,4.9765369705661607e-17,3.3881317890172006e-21,7.79002343368992e-06,1.6364022961483934e-17,0.06446054173482972,7.79002343368992e-06
1,100.0,3.7317930012595747e-17,1.64464618860614e-18,3.866535256984768,0.0,0.00010000000000000018,optimaltree,10.0,1.0638049822737197e-16,1.64464618860614e-16,12.428490091836007,0.02068065166473388,0.00013559103012084962,0.08048673629760743,0.020816242694854728,20.0,1.0,270.0,100.0,0.0,1.0,100.0,10.0,10000.0,220.0,OCT,facility,4.9765369705661607e-17,1.0408340855860843e-17,6.465083237246323e-06,1.6364022961483934e-17,0.061305855962301335,6.465083237246323e-06
2,100.0,3.7317930012595747e-17,1.64464618860614e-18,573.2442040620863,0.0,0.00010000000000000018,optimaltree,10.0,1.0638049822737197e-16,1.64464618860614e-16,1912.755484339516,2.8007030487060547e-05,0.00014126539230346682,0.0970344352722168,0.00016927242279052738,20.0,1.0,270.0,100.0,0.0,1.0,100.0,10.0,10000.0,220.0,OCT-H,facility,4.9765369705661607e-17,3.3881317890172006e-21,8.74255247339638e-06,1.6364022961483934e-17,0.0826778997543652,8.74255247339638e-06
0,100.0,9.451083287828508e-07,-2.5968783198022643e-09,149.05803453187974,0.0001,0.00020098000000000018,pytorch,20.0,9.451083286723826e-05,1.81344048669619e-09,310.0417272705596,1.0931491851806648e-05,0.0013213229179382324,0.19858322381973267,0.0013322544097900387,20.0,9.0,480.0,100.0,0.0,9.0,100.0,20.0,10000.0,420.0,NN,facility,9.403709137327447e-06,6.776263578034403e-21,4.1002609417627546e-05,2.6400867922514035e-08,0.07448198166956534,4.1002609417627546e-05
1,100.0,9.451083287828508e-07,-2.5968783198022643e-09,174.62197216763522,0.0001,0.00020098000000000018,optimaltree,20.0,9.451083286723826e-05,1.81344048669619e-09,298.9262088339757,4.790544509887696e-05,0.0013206362724304196,0.2389774537086487,0.0013685417175292975,20.0,9.0,480.0,100.0,0.0,9.0,100.0,20.0,10000.0,420.0,OCT,facility,9.403709137327447e-06,6.776263578034403e-21,3.863827962686085e-05,2.6400867922514035e-08,0.08943894732197938,3.863827962686085e-05
2,100.0,9.451083287828508e-07,-2.5968783198022643e-09,234.11336949740158,0.0001,0.00020098000000000018,optimaltree,20.0,9.451083286723826e-05,1.81344048669619e-09,1068.0004643333914,4.7957897186279325e-05,0.0012283873558044433,0.29880948781967165,0.0012763452529907231,20.0,9.0,480.0,100.0,0.0,9.0,100.0,20.0,10000.0,420.0,OCT-H,facility,9.403709137327447e-06,2.7105054312137605e-20,4.70124060279583e-05,2.6400867922514035e-08,0.3190413540558095,4.70124060279583e-05
0,100.0,1.0487887533349931e-16,-1.2145010640699696e-17,1986.001796282022,0.0,0.00010000000000000018,pytorch,20.0,1.7639640319923504e-16,2.4261343556639216e-16,2618.013145852401,1.0261535644531253e-05,0.0003547430038452149,0.7248996710777282,0.0003650045394897462,40.0,2.0,940.0,100.0,0.0,2.0,100.0,20.0,10000.0,840.0,NN,facility,6.9535665773389e-17,1.6940658945086007e-21,1.565705900267701e-05,1.0497923399109683e-16,0.2278209780355729,1.565705900267701e-05
1,100.0,1.0487887533349931e-16,-1.2145010640699696e-17,2080.121600988587,0.0,0.00010000000000000018,optimaltree,20.0,1.7639640319923504e-16,2.4261343556639216e-16,5646.144200950469,3.909826278686524e-05,0.00036221027374267584,0.8347705554962158,0.000401308536529541,40.0,2.0,940.0,100.0,0.0,2.0,100.0,20.0,10000.0,840.0,OCT,facility,6.9535665773389e-17,0.0,1.2000608231914458e-05,1.0497923399109683e-16,0.3765629204156951,1.2000608231914458e-05
2,100.0,1.0487887533349931e-16,-1.2145010640699696e-17,2201.2310199309213,0.0,0.00010000000000000018,optimaltree,20.0,1.7639640319923504e-16,2.4261343556639216e-16,3991.7192130300127,3.939628601074219e-05,0.00031540393829345705,0.7809972596168518,0.0003548002243041992,40.0,2.0,940.0,100.0,0.0,2.0,100.0,20.0,10000.0,840.0,OCT-H,facility,6.9535665773389e-17,0.0,1.527619486718529e-05,1.0497923399109683e-16,0.26601774780461124,1.527619486718529e-05
0,100.0,9.873971094363796e-17,6.581341424441165e-11,1037.0305152028054,0.0001,0.00020296000000000016,pytorch,40.0,1.4183724604570051e-16,2.5014226539824282e-09,1522.7221429413073,9.787082672119143e-06,0.0007733988761901855,0.812187738418579,0.0007831859588623048,40.0,5.0,1760.0,100.0,0.0,5.0,100.0,40.0,10000.0,1640.0,NN,facility,2.0088399620467625e-17,3.3881317890172006e-21,4.2252678282216e-05,3.0670686433659506e-10,0.2419722803184555,4.2252678282216e-05
1,100.0,9.873971094363796e-17,6.581341424441165e-11,961.9083487293512,0.0001,0.00020296000000000016,optimaltree,40.0,1.4183724604570051e-16,2.5014226539824282e-09,1392.698164952522,4.1606426239013656e-05,0.0008333373069763184,0.8416156816482544,0.0008749437332153324,40.0,5.0,1760.0,100.0,0.0,5.0,100.0,40.0,10000.0,1640.0,OCT,facility,2.0088399620467625e-17,1.3552527156068802e-20,6.325460002353376e-05,3.0670686433659506e-10,0.25194607630150345,6.325460002353376e-05
2,100.0,9.873971094363796e-17,6.581341424441165e-11,1070.2035796932885,0.0001,0.00020296000000000016,optimaltree,40.0,1.4183724604570051e-16,2.5014226539824282e-09,1680.8184784091895,4.0924549102783225e-05,0.0007544493675231932,0.8512120127677918,0.000795373916625977,40.0,5.0,1760.0,100.0,0.0,5.0,100.0,40.0,10000.0,1640.0,OCT-H,facility,2.0088399620467625e-17,2.7105054312137605e-20,3.5034925032928556e-05,3.0670686433659506e-10,0.2689088279428047,3.5034925032928556e-05
0,100.0,9.874335541819171e-17,1.3240439232523864e-17,325.76148189385026,0.0001,0.00020098000000000018,pytorch,30.0,2.0202332213011475e-16,1.9051789454449635e-16,519.1451683428072,1.082420349121093e-05,0.0005202341079711914,0.1729983425140381,0.0005310583114624024,60.0,3.0,2010.0,100.0,0.0,3.0,100.0,30.0,10000.0,1860.0,NN,facility,7.883798368082264e-17,6.776263578034403e-21,2.986411326508908e-05,1.0783246977224071e-16,0.05725944035889061,2.986411326508908e-05
1,100.0,9.874335541819171e-17,1.3240439232523864e-17,334.57772090993774,0.0001,0.00020098000000000018,optimaltree,30.0,2.0202332213011475e-16,1.9051789454449635e-16,912.92994731485,3.724575042724609e-05,0.0005692696571350098,0.20292654275894165,0.0006065154075622559,60.0,3.0,2010.0,100.0,0.0,3.0,100.0,30.0,10000.0,1860.0,OCT,facility,7.883798368082264e-17,6.776263578034403e-21,4.1253920172583007e-05,1.0783246977224071e-16,0.09347167136148052,4.1253920172583e-05
2,100.0,9.874335541819171e-17,1.3240439232523864e-17,296.35777599017916,0.0001,0.00020098000000000018,optimaltree,30.0,2.0202332213011475e-16,1.9051789454449635e-16,653.6152935518478,4.04644012451172e-05,0.0005849099159240723,0.18533454179763795,0.0006253743171691895,60.0,3.0,2010.0,100.0,0.0,3.0,100.0,30.0,10000.0,1860.0,OCT-H,facility,7.883798368082264e-17,6.776263578034403e-21,3.107702351424717e-05,1.0783246977224071e-16,0.06984379940584894,3.107702351424717e-05
0,100.0,1.8204712684545407e-16,3.931949187107195e-18,341.2440108518992,0.0001,0.00019900000000000018,pytorch,60.0,3.0138522097093047e-16,2.62821543995524e-16,552.638892701271,1.0719299316406247e-05,0.0011290955543518067,0.3889549922943115,0.001139814853668213,60.0,5.0,3840.0,100.0,0.0,5.0,100.0,60.0,10000.0,3660.0,NN,facility,5.1259177963570235e-17,1.6940658945086007e-21,4.1675088263729276e-05,1.0064186771015923e-16,0.1320573999089376,4.1675088263729276e-05
1,100.0,1.8204712684545407e-16,3.931949187107195e-18,415.5359579598195,0.0001,0.00019900000000000018,optimaltree,60.0,3.0138522097093047e-16,2.62821543995524e-16,1576.854420120437,4.134893417358399e-05,0.0011432528495788571,0.4922446370124817,0.0011846017837524414,60.0,5.0,3840.0,100.0,0.0,5.0,100.0,60.0,10000.0,3660.0,OCT,facility,5.1259177963570235e-17,6.776263578034403e-21,4.9213249129531876e-05,1.0064186771015923e-16,0.4463303890350442,4.9213249129531876e-05
2,100.0,1.8204712684545407e-16,3.931949187107195e-18,409.0837105759635,0.0001,0.00019900000000000018,optimaltree,60.0,3.0138522097093047e-16,2.62821543995524e-16,1447.4046923695093,4.9555301666259744e-05,0.0014530324935913086,0.614684190750122,0.001502587795257568,60.0,5.0,3840.0,100.0,0.0,5.0,100.0,60.0,10000.0,3660.0,OCT-H,facility,5.1259177963570235e-17,2.7105054312137605e-20,4.716026070627694e-05,1.0064186771015923e-16,0.5008000868795766,4.716026070627695e-05
0,100.0,1.6654094778714536e-16,1.7977455950049813e-17,97.37588950857331,0.0001,0.00020098000000000018,pytorch,40.0,3.0250164518631123e-16,2.770987445478483e-16,176.02862196956283,1.1830329895019529e-05,0.002572352886199951,0.2516371393203735,0.0025841832160949712,80.0,9.0,3480.0,100.0,0.0,9.0,100.0,40.0,10000.0,3280.0,NN,facility,4.992267735817712e-17,1.6940658945086007e-21,9.181539942930342e-05,1.0802067325567659e-16,0.07114337563967398,9.181539942930342e-05
1,100.0,1.6654094778714536e-16,1.7977455950049813e-17,99.44346492037752,0.0001,0.00020098000000000018,optimaltree,40.0,3.0250164518631123e-16,2.770987445478483e-16,166.73776157487262,5.4388046264648405e-05,0.0025673270225524897,0.2607124304771423,0.002621715068817139,80.0,9.0,3480.0,100.0,0.0,9.0,100.0,40.0,10000.0,3280.0,OCT,facility,4.992267735817712e-17,2.7105054312137605e-20,8.803275598666633e-05,1.0802067325567659e-16,0.0783874399533762,8.803275598666633e-05
2,100.0,1.6654094778714536e-16,1.7977455950049813e-17,89.51884012939243,0.0001,0.00020098000000000018,optimaltree,40.0,3.0250164518631123e-16,2.770987445478483e-16,191.8504803961348,4.1878223419189475e-05,0.0025362896919250487,0.2307946014404297,0.002578167915344237,80.0,9.0,3480.0,100.0,0.0,9.0,100.0,40.0,10000.0,3280.0,OCT-H,facility,4.992267735817712e-17,2.7105054312137605e-20,7.827027737506475e-05,1.0802067325567659e-16,0.0884311500249438,7.827027737506475e-05
0,100.0,1.900292926805707e-16,2.5171949757471797e-12,865.4516848761031,0.0001,0.00020098000000000018,pytorch,80.0,3.3396972375238826e-16,2.517189640599084e-10,1524.125299174308,1.605749130249024e-05,0.001444699764251709,1.2642148280143737,0.001460757255554199,80.0,6.0,6720.0,100.0,0.0,6.0,100.0,80.0,10000.0,6480.0,NN,facility,5.5099467207339915e-17,6.776263578034403e-21,0.0001717745669201785,2.504572015605215e-11,0.46224824904841094,0.00017177456692017847
1,100.0,1.900292926805707e-16,2.5171949757471797e-12,64.73933922290652,0.0001,0.00020098000000000018,optimaltree,80.0,3.3396972375238826e-16,2.517189640599084e-10,199.64562979966993,0.0196616005897522,0.0013933515548706055,1.363083689212799,0.02105495214462281,80.0,6.0,6720.0,100.0,0.0,6.0,100.0,80.0,10000.0,6480.0,OCT,facility,5.5099467207339915e-17,3.469446951953614e-18,0.00017760651803964636,2.504572015605215e-11,0.7459456041796121,0.00017760651803964633
2,100.0,1.900292926805707e-16,2.5171949757471797e-12,1061.2155014854047,0.0001,0.00020098000000000018,optimaltree,80.0,3.3396972375238826e-16,2.517189640599084e-10,2971.5480481384147,7.258892059326171e-05,0.0014040803909301759,1.5670643639564514,0.001476669311523437,80.0,6.0,6720.0,100.0,0.0,6.0,100.0,80.0,10000.0,6480.0,OCT-H,facility,5.5099467207339915e-17,0.0,0.00014882648419869867,2.504572015605215e-11,1.2924675952507079,0.00014882648419869867
\end{filecontents*}
\csvreader[head to column names, late after line=\\]{./data/benchmarks/facility/facility_general.csv}{
predictor=\predictor,
accuracy=\accuracy,
n_strategies=\nstrategies,
n_infeas=\ninfeas,
avg_infeas=\avginfeas,
avg_subopt=\avgsubopt,
max_infeas=\maxinfeas,
max_subopt=\maxsubopt,
avg_time_improv=\avgtimeimprov,
max_time_improv=\maxtimeimprov,
n=\n,
m=\m,
n_var=\nvar,
n_constr=\nconstr,
n_train=\ntrain,
good_turing=\goodturing,
}
{\tablenum[table-format=2, round-precision=0]{\n} & \m & \nvar & \nconstr & \predictor & \ntrain & \goodturing &\nstrategies &\accuracy & \avgtimeimprov & \maxinfeas & \maxsubopt } 
     \bottomrule
\end{tabular}
\end{adjustbox}
\end{table}

\subsection{Hybrid Vehicle Control}%
\label{ssub:hybrid_vehicle_control}
Consider the hybrid-vehicle control problem in~\cite[\Sec 3.2]{takapoui2017}.
The model consists of a battery, an electric motor/generator and an engine.
We assume to know the demanded power $P^{\rm des}_t$ over the horizon $t=0,\dots,T-1$.
The goal is to plan the battery and the engine power outputs $P_{t}^{\rm batt}$ and $P^{\rm eng}_{t}$ so that they match at least the demand,
\begin{equation*}
    P_{t}^{\rm batt} + P_{t}^{\rm eng} \ge P_{t}^{\rm des}.
\end{equation*}
The status of the battery is modeled with the internal energy $E_t$ evolving as
\begin{equation*}
    E_{t+1} = E_{t} - \tau P_{t}^{\rm batt},
\end{equation*}
where $\tau > 0$ is the time interval discretization.
The battery capacity is limited by $E^{\rm max}$ and its initial value is $E_{\rm init}$.

We now model the cost function.
We penalize the terminal energy state of the battery with the function
\begin{equation*}
    g(E) = \eta (E^{\rm max} - E)^2.
\end{equation*}
with $\eta \ge 0$.
At each time $t$ we model the on-off state of the engine with the binary variable $z_t$.
When the engine is off ($z_t = 0$) we do not consume any energy, thus we have $0 \le P_{t}^{\rm eng} \le P^{\rm max} z_t$.
When the engine is on ($z_t = 1$) it consumes $\alpha P_{t}^{\rm eng} + \beta P_{t}^{\rm eng} + \gamma$ units of fuel, with $\alpha, \beta, \gamma > 0$.
We define the stage power cost as
\begin{equation*}
    f(P, z) = \alpha P^2 + \beta P + \gamma z.
\end{equation*}
We also introduce a cost of turning the engine on at time $t$ as $\delta (z_t - z_{t-1})_{+}$.

The hybrid vehicle control problem can be formulated as a mixed integer quadratic optimization problem:
\begin{equation*}
\begin{array}{ll}
    \mbox{minimize}   &  \displaystyle \eta (E_T - E^{\rm max})^2 + \sum_{t=0}^{T - 1} f(P_{t}^{\rm eng}, z_t) + \delta (z_t - z_{t-1})_{+}\\
    \mbox{subject to} & E_{t + 1} = E_{t} - \tau P_{t}^{\rm batt}, \qquad t=0,\dots,T-1\\
    & 0 \le E_t \le E^{\rm max}, \qquad t=0,\dots,T\\
    & E_0 = E_{\rm init} \\
    & 0 \le P^{\rm eng}_t \le P^{\rm max}, \qquad t=0, \dots, T-1\\
    & P_{t}^{\rm batt} + P_{t}^{\rm eng} \ge P_{t}^{\rm des}, \qquad t=0,\dots,T-1\\
    & z_t \in \{0, 1\}, \qquad t=0,\dots,T-1.
\end{array}
\end{equation*}
The initial state $E_0$ and demanded power $P^{\rm des}_t$ are our parameters.

\paragraph{Problem instances.}
We generated the control instances with varying horizon length $T$.
The discretization time interval is $\tau = 4$.
We chose the cost function parameters $\alpha=\beta=\gamma=1$, and $\delta = 0.1$.
The maximum charge is $E^{\rm max}= 50$ and the maximum power $P^{\rm max}=1$.
The parameter $E_0$ was sampled from a uniform distribution within the  ball $B(40, 0.5)$.
The demand $P^{\rm des}$ also comes from a uniform  distribution within the ball  $B(\overline{P}^{\rm des}, 0.5)$ where
\ifpreprint
\begin{align*}
\overline{P}^{\rm des}_t =\; &(0.05, 0.30, 0.55, 0.80, 1.05, 1.30, 1.55, 1.80, 1.95, 1.70, 1.45, 1.20, 1.02, \\
&1.12, 1.22, 1.32, 1.42, 1.52, 1.62, 1.72, 1.73, 1.38, 1.03, 0.68, 0.33, -0.02, \\
&-0.37, -0.72, -0.94, -0.64, -0.34, -0.04, 0.18, 0.08, -0.02, -0.12,\\
&-0.22, -0.32, -0.42, -0.52)
\end{align*}
\fi
is the desired power in~\cite{takapoui2017}.
Depending on the horizon length $T$ we choose only the first $T$ elements of $\overline{P}^{\rm des}$ to sample $P^{\rm des}$.


\paragraph{Results.}
We present  the results in Table~\ref{tab:control}.
Here there is a significant improvement in computation time compared to Gurobi.
Also in this case the number of strategies $|\strategy|$ is much less than the number of variables or the worst case number of control input combinations.
The maximum infeasibility and suboptimality are  low even when the predictions are not exact and the performance of \glspl{OCT}, \glspl{OCT-H} and NNs is comparable except in the last case.


\begin{table}
  \centering
    \caption{Hybrid control benchmarks.}
    \label{tab:control}
    \begin{adjustbox}{max width=\textwidth}
\begin{tabular}{l
                S[table-format=3, round-precision=0]
                S[table-format=3, round-precision=0]
        l
                S[table-format=5, round-precision=0]
                S[table-format = 1.2e1,
                  scientific-notation=true,
                  retain-zero-exponent=true,
                  round-integer-to-decimal]
                S[table-format=2, round-precision=0]
                S[table-format=3.2, round-integer-to-decimal]
                S[table-format=4.2, round-integer-to-decimal]
                S[table-format = 1.2e2,
                  scientific-notation=true,
                  retain-zero-exponent=true,
                  round-integer-to-decimal]
                S[table-format = 1.2e2,
                  scientific-notation=true,
                  retain-zero-exponent=true,
                  round-integer-to-decimal]
                }
\toprule
$T$ & $n_{\rm var}$ & $n_{\rm con}$ & learner & $N$ & ${\rm GT}$ & $|\strategy|$ & ${\rm acc}\;[\%]$ & $t_{\rm ratio}$ & ${\bar{p}}$ & ${\bar{d}}$ \\
\midrule
\begin{filecontents*}{./data/benchmarks/control/control_general.csv}
,T,accuracy,avg_infeas,avg_subopt,avg_time_improv,good_turing,good_turing_smooth,learner,max_infeas,max_subopt,max_time_improv,mean_pred_time_pred,mean_solve_time_pred,mean_time_full,mean_time_pred,n_best,n_constr,n_correct,n_infeas,n_strategies,n_test,n_theta,n_train,n_var,predictor,problem,std_infeas,std_pred_time_pred,std_solve_time_pred,std_subopt,std_time_full,std_time_pred
0,10.0,100.0,5.79471665914379e-16,1.901261942216864e-15,2.790813688999173,0.0003,0.0004029400000000002,pytorch,1.1689028438955087e-14,3.777265400106323e-15,3.579204827276455,1.3761520385742194e-05,0.0008163928985595703,0.0023168063163757323,0.0008301544189453124,10.0,103.0,100.0,0.0,16.0,100.0,11.0,10000.0,51.0,NN,control,1.5198725925036318e-15,6.776263578034403e-21,3.794789522108281e-05,7.90573652463569e-16,0.0004792939907909172,3.794789522108282e-05
1,10.0,99.0,9.365109757471255e-06,1.8992595218823212e-15,0.12358584602645212,0.0003,0.0004029400000000002,optimaltree,0.0009365109756893283,3.777265400106323e-15,0.8703434509458537,0.02322510480880738,0.0008107805252075195,0.0029704952239990234,0.024035885334014895,10.0,103.0,99.0,1.0,16.0,100.0,11.0,10000.0,51.0,OCT,control,9.318166555227639e-05,6.938893903907228e-18,4.593474364120456e-05,7.943040540243771e-16,0.0022877484239083477,4.593474364120457e-05
2,10.0,100.0,5.79471665914379e-16,1.901261942216864e-15,3.0024586222408365,0.0003,0.0004029400000000002,optimaltree,1.1689028438955087e-14,3.777265400106323e-15,11.477931197022276,0.00010559082031249999,0.0008166170120239258,0.002768890857696533,0.0009222078323364258,10.0,103.0,100.0,0.0,16.0,100.0,11.0,10000.0,51.0,OCT-H,control,1.5198725925036318e-15,1.3552527156068805e-20,5.1875500152048785e-05,7.90573652463569e-16,0.0015412154295582662,5.1875500152048785e-05
0,15.0,99.0,7.267116673646766e-07,-1.3709469847881219e-06,4.204393974848753,0.0004,0.0005039200000000001,pytorch,7.267116658959505e-05,4.2871033554385475e-15,15.056212138941644,1.1177062988281252e-05,0.0008956217765808106,0.003812539577484131,0.0009067988395690918,10.0,153.0,99.0,0.0,16.0,100.0,16.0,10000.0,76.0,NN,control,7.2306897795053275e-06,1.6940658945086007e-21,4.038102019944848e-05,1.3640750293129283e-05,0.002283532257269274,4.038102019944848e-05
1,15.0,99.0,7.267116673646766e-07,-1.3709469847881219e-06,3.675778533436289,0.0004,0.0005039200000000001,optimaltree,7.267116658959505e-05,4.2871033554385475e-15,15.93720144122221,9.06610488891602e-05,0.0009035491943359376,0.003654496669769287,0.000994210243225098,10.0,153.0,99.0,0.0,16.0,100.0,16.0,10000.0,76.0,OCT,control,7.2306897795053275e-06,4.0657581468206416e-20,4.0051328456420245e-05,1.3640750293129283e-05,0.0015462437821545602,4.0051328456420245e-05
2,15.0,99.0,7.267116673646766e-07,-1.3709469847881219e-06,3.0522233426675545,0.0004,0.0005039200000000001,optimaltree,7.267116658959505e-05,4.2871033554385475e-15,8.790118618223238,8.723258972167967e-05,0.0011634945869445801,0.0038174986839294434,0.0012507271766662594,10.0,153.0,99.0,0.0,16.0,100.0,16.0,10000.0,76.0,OCT-H,control,7.2306897795053275e-06,1.3552527156068805e-20,0.0001729929406552371,1.3640750293129283e-05,0.0014340346083045066,0.0001729929406552371
0,20.0,100.0,8.97385622034149e-16,-1.7588746373718847e-07,3.888838834951456,0.0006,0.0006999400000000001,pytorch,1.4636107759479117e-14,6.043923096149303e-15,4.846001820886905,1.4853477478027343e-05,0.0009060382843017578,0.0035811996459960936,0.0009208917617797852,10.0,203.0,100.0,0.0,16.0,100.0,21.0,10000.0,101.0,NN,control,2.0439097655529143e-15,0.0,5.0442509013198937e-05,1.7500582081108618e-06,0.0006031010820794897,5.0442509013198937e-05
1,20.0,99.0,2.2294020001876494e-06,-1.7766410481747587e-07,4.15572900000464,0.0006,0.0006999400000000001,optimaltree,0.00022294019992936725,6.043923096149303e-15,12.432543720379616,9.809017181396482e-05,0.0009301185607910156,0.004272956848144532,0.00102820873260498,10.0,203.0,99.0,1.0,16.0,100.0,21.0,10000.0,101.0,OCT,control,2.218226981548601e-05,2.710505431213761e-20,4.8996764056796744e-05,1.758784946420243e-06,0.001704156047055094,4.899676405679673e-05
2,20.0,100.0,8.97385622034149e-16,-1.7588746373718847e-07,4.333670467519056,0.0006,0.0006999400000000001,optimaltree,1.4636107759479117e-14,6.043923096149303e-15,17.87458697065108,0.00010200977325439453,0.0009045600891113281,0.004362142086029053,0.0010065698623657226,10.0,203.0,100.0,0.0,16.0,100.0,21.0,10000.0,101.0,OCT-H,control,2.0439097655529143e-15,0.0,4.5656516715545003e-05,1.7500582081108618e-06,0.002467679432951981,4.5656516715544976e-05
0,25.0,100.0,7.792759827731085e-07,-1.3922835337729376e-06,4.2672990923445955,0.0002,0.0003019600000000002,pytorch,5.489883775435924e-05,6.326408474558817e-15,5.387969780208294,1.3976097106933591e-05,0.0009749960899353027,0.004220240116119385,0.0009889721870422364,10.0,253.0,100.0,0.0,10.0,100.0,26.0,10000.0,126.0,NN,control,5.902100549150255e-06,1.6940658945086007e-21,5.040852454982971e-05,9.463316304154256e-06,0.0006015038725133298,5.040852454982971e-05
1,25.0,100.0,7.792759827731085e-07,-1.3922835337729376e-06,4.3787736910818,0.0002,0.0003019600000000002,optimaltree,5.489883775435924e-05,6.326408474558817e-15,6.99390567566381,7.825374603271488e-05,0.0009809446334838867,0.004637989997863769,0.0010591983795166014,10.0,253.0,100.0,0.0,10.0,100.0,26.0,10000.0,126.0,OCT,control,5.902100549150255e-06,2.710505431213761e-20,4.246096043830603e-05,9.463316304154256e-06,0.0006798402618248867,4.246096043830603e-05
2,25.0,100.0,7.792759827731085e-07,-1.3922835337729376e-06,5.119721396255165,0.0002,0.0003019600000000002,optimaltree,5.489883775435924e-05,6.326408474558817e-15,13.69254943040461,0.00011441707611083984,0.0009463858604431152,0.005431015491485596,0.0010608029365539555,10.0,253.0,100.0,0.0,10.0,100.0,26.0,10000.0,126.0,OCT-H,control,5.902100549150255e-06,0.0,5.386347193848697e-05,9.463316304154256e-06,0.0030288603284382466,5.386347193848696e-05
0,30.0,100.0,2.7956240612103196e-15,-1.9297012968486793e-16,86.6532050429429,0.0001,0.00020098000000000018,pytorch,1.548196564810179e-13,7.720885664475351e-16,161.33434224271818,1.2192726135253908e-05,0.0009916067123413087,0.08698243856430053,0.0010037994384765627,10.0,303.0,100.0,0.0,20.0,100.0,31.0,10000.0,151.0,NN,control,1.5858410437979327e-14,1.6940658945086007e-21,6.0897449215583896e-05,3.7211489994652297e-16,0.040070756585911985,6.089744921558393e-05
1,30.0,100.0,2.793731417458717e-15,-1.8636633493282847e-16,148.36844187667197,0.0001,0.00020098000000000018,optimaltree,1.548196564810179e-13,7.720885664475351e-16,774.8570726675975,8.679151535034176e-05,0.0009737658500671387,0.15735324382781982,0.0010605573654174805,10.0,303.0,100.0,0.0,20.0,100.0,31.0,10000.0,151.0,OCT,control,1.5858710674083237e-14,2.710505431213761e-20,5.887873023210228e-05,3.7108638984444313e-16,0.20734251459280742,5.887873023210228e-05
2,30.0,100.0,2.7885822572161943e-15,-1.913173947248919e-16,124.01485672489453,0.0001,0.00020098000000000018,optimaltree,1.548196564810179e-13,7.720885664475351e-16,680.7726686998388,0.0001190352439880371,0.000989229679107666,0.13744131565093995,0.0011082649230957037,10.0,303.0,100.0,0.0,20.0,100.0,31.0,10000.0,151.0,OCT-H,control,1.5859495064235458e-14,0.0,5.22720286744164e-05,3.711318603468611e-16,0.16245045908526715,5.22720286744164e-05
0,35.0,100.0,4.0413565961941343e-07,-2.511198365550096e-10,845.128667063275,0.0003,0.00040888000000000014,pytorch,4.041356587741283e-05,1.136672018044222e-15,1151.2127003881826,1.0852813720703127e-05,0.001067197322845459,0.9110910749435425,0.0010780501365661626,10.0,353.0,100.0,0.0,37.0,100.0,36.0,10000.0,176.0,NN,control,4.021099033585699e-06,1.6940658945086007e-21,5.884281987992558e-05,2.498607563250245e-09,0.3071099466213858,5.884281987992557e-05
1,35.0,98.0,4.041356594614654e-07,6.342531303628582e-06,1010.7264127926571,0.0003,0.00040888000000000014,optimaltree,4.041356587741283e-05,0.0003479718941274592,3971.8592726959205,8.613348007202145e-05,0.0011940145492553712,1.2938794255256654,0.0012801480293273932,10.0,353.0,98.0,0.0,37.0,100.0,36.0,10000.0,176.0,OCT,control,4.0210990336015735e-06,2.710505431213761e-20,5.126973902321052e-05,3.9025218591699964e-05,0.9574892970258215,5.126973902321052e-05
2,35.0,100.0,4.041356595075794e-07,5.977364304028585e-07,1250.1897275514964,0.0003,0.00040888000000000014,optimaltree,4.041356587741283e-05,5.979875502232859e-05,4760.884165712419,7.738590240478517e-05,0.00105621337890625,1.4172141766548156,0.0011335992813110356,10.0,353.0,100.0,0.0,37.0,100.0,36.0,10000.0,176.0,OCT-H,control,4.021099033596939e-06,1.3552527156068805e-20,5.776435604145258e-05,5.949926763253046e-06,1.34354997885749,5.776435604145258e-05
0,40.0,100.0,1.4679614879603982e-15,-2.079534651816084e-16,6461.459786523925,0.0002,0.00029998000000000015,pytorch,6.498071577878773e-14,1.5316488156865927e-15,9526.168392957003,1.8699169158935555e-05,0.0009549617767333985,6.291271047592163,0.000973660945892334,10.0,403.0,100.0,0.0,35.0,100.0,41.0,10000.0,201.0,NN,control,6.79008099276823e-15,6.776263578034403e-21,3.8091141852567085e-05,4.3493166692624977e-16,1.5186156884238737,3.809114185256713e-05
1,40.0,93.0,1.41529491065667e-15,0.000637205004633362,6846.486822515289,0.0002,0.00029998000000000015,optimaltree,6.498071577878773e-14,0.031199629474100266,22038.071971781093,0.00015332221984863282,0.0009740447998046875,7.718503444194794,0.0011273670196533206,10.0,403.0,93.0,0.0,35.0,100.0,41.0,10000.0,201.0,OCT,control,6.790707964585691e-15,0.0,4.492329052727274e-05,0.00429840593645891,4.915845492348019,4.492329052727274e-05
2,40.0,89.0,1.287657499700929e-15,0.0003853845764890452,6551.657296646687,0.0002,0.00029998000000000015,optimaltree,6.498071577878773e-14,0.030228625184043786,13400.799599507696,8.506298065185544e-05,0.0009766650199890137,6.95607800245285,0.0010617280006408694,10.0,403.0,89.0,0.0,35.0,100.0,41.0,10000.0,201.0,OCT-H,control,6.482771632236728e-15,2.710505431213761e-20,4.625709083539525e-05,0.0030173198383965606,2.5680240300602994,4.6257090835395244e-05
\end{filecontents*}
\csvreader[head to column names, late after line=\\]{./data/benchmarks/control/control_general.csv}{
predictor=\predictor,
accuracy=\accuracy,
n_strategies=\nstrategies,
n_infeas=\ninfeas,
avg_infeas=\avginfeas,
avg_subopt=\avgsubopt,
max_infeas=\maxinfeas,
max_subopt=\maxsubopt,
avg_time_improv=\avgtimeimprov,
max_time_improv=\maxtimeimprov,
T=\T,
n_var=\nvar,
n_constr=\nconstr,
n_train=\ntrain,
good_turing=\goodturing,
}
{\tablenum[table-format=2, round-precision=0]{\T} & \nvar & \nconstr & \predictor & \ntrain & \goodturing &\nstrategies &\accuracy & \avgtimeimprov & \maxinfeas & \maxsubopt } 
     \bottomrule
\end{tabular}
\end{adjustbox}
\end{table}

\section{Conclusions}
\label{sec:conclusions}

We introduced the idea that using \glspl{OCT} and \glspl{OCT-H} we  obtain insight on the strategy of the optimal solution  in  many optimization problems
as a function of key parameters.
In this way,  optimization is not a black box anymore, but rather it has a voice, \ie,  we are able to provide insights on the logic behind the optimal solution. The class of optimization problems that these ideas apply to is rather broad since it includes any continuous and mixed-integer convex optimization problems without any assumption on the parameters dependency.
The accuracy of our approach is in the 90\%-100\%, while even when it does not provide the optimal solution the degree of suboptimality or infeasibility is extremely low.
In addition, the computation times of our method are significantly faster than solving the problem using standard optimization algorithms.
Even though the solving times did not have a strict limit for the proposed applications, our method can be useful in online settings where problems must be solved within a limited time and with limited computational resources. Applications of this framework to fast online optimization appear in~\citep{bertsimas2019}.
Comparisons on several examples show that the out-of-sample strategy prediction accuracy of \glspl{OCT-H} is comparable to the \glspl{NN}.
\glspl{OCT} show competitive accuracy in most comparisons apart from the cases with higher number of strategies where the prediction performance degrades.
Note that this is not necessarily related to the size of the problem but rather to how the parameters affect the optimal solution and, in turn, the total number of strategies.
\glspl{OCT} also exhibit low infeasibility and low suboptimality even when the prediction is not correct.
In addition, \glspl{OCT} have higher interpretability than~\glspl{OCT-H} and \glspl{NN} because of their structure with parallel splits, therefore, being required in applications where models must be explainable.
For these reasons, our framework provides a reliable and insightful understanding of optimal solutions in a broad class of optimization problems using different machine learning algorithms.



\appendix

\section{Optimal  Classification Trees} \label{sec:trees}
\glspl{OCT} and \glspl{OCT-H} developed by~\citet{bertsimas2017a,bd_book} are a recently proposed  generalization of~\glspl{CART} developed by~\citet{breiman1984} that
construct decision trees that are near  optimal with significantly higher prediction accuracy while retaining their interpretability.
\citet{bms} have shown that  given a \gls{NN} (feedforward, convolutional and recurrent), we can construct an \gls{OCT-H} that has the same in sample accuracy, that is \glspl{OCT-H} are at least as powerful as \glspl{NN}. The constructions can sometimes generate OCT-Hs with high depth.
However, they also report computational results that show that \glspl{OCT-H} and \glspl{NN} have comparable performance in practice even for depths of \glspl{OCT-H} below 10.

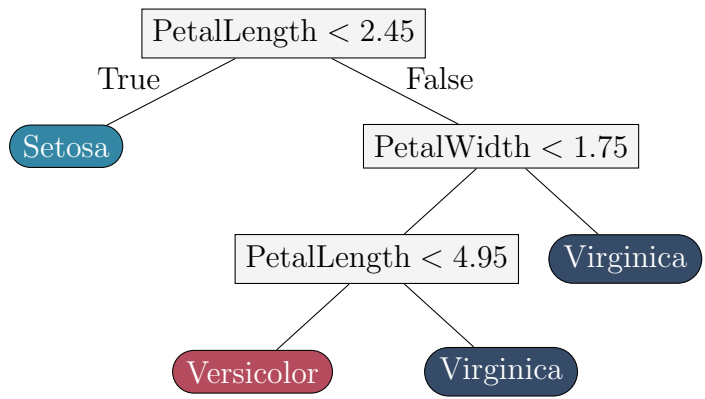
\begin{figure}
    \centering
    \begin{tikzpicture}[sibling distance=10em, level 3/.style={sibling distance=8em}]
    \node[split] {${\rm Petal Length} < 2.45$}
        child {node[class3]{Setosa}
               edge from parent node[left, yshift=1ex] {True}}
        child { node[split] {${\rm Petal Width} < 1.75$}
                    child {node[split] {${\rm Petal Length} < 4.95$}
                        child {node[class1]{Versicolor}}
                        child {node[class2]{Virginica}}
                    }
                    child {node[class2]{Virginica}}
                edge from parent node[right, yshift=1ex] {False}}
    ;
\end{tikzpicture}
    \caption{Example \gls{OCT} for the Iris dataset~\cite{bd_book}.}
    \label{fig:oct}
\end{figure}

\paragraph{Architecture.}
\Glspl{OCT} recursively partition the feature space $\reals^p$ to construct hierarchical disjoint regions.
A tree can be defined as a set of nodes $t \in \tree$ of two types $\tree = \tree_B \cup \tree_L$:
\begin{description}
\item [Branch nodes] Nodes $t \in \tree_B$ at the tree branches describe a split of the form ${a_t^\tpose \theta < b_t}$ where $a_t \in \reals^p$ and $b_t \in \reals$. They partition the space in two subsets: the points on the left branch satisfying the inequality and the remaining ones points on the right branch.
If splits involve a single variable we denote them as \emph{parallel} and we refer to the tree as \glsfirst{OCT}.
This is achieved by enforcing all  components of $a_t$ to be all $0$ except from one.
Otherwise, if the components of $a_t$ can be freely nonzero, we denote the splits as \emph{hyperplanes} and we refer to the tree as \glsfirst{OCT-H}.
\item [Leaf nodes] Nodes $t \in \tree_L$ at the tree leaves make a class prediction for each point falling into that node.
\end{description}
An example \gls{OCT} for the Iris dataset appears in Figure~\ref{fig:oct}~\citep{bd_book}.
For each new data point it is straightforward to follow which hyperplanes are satisfied and to make a prediction.
This characteristic makes \glspl{OCT} and \glspl{OCT-H} highly interpretable.
Note that the level of interpretability of the resulting trees can be also tuned by changing minimum sparsity of $a_t$. The two extremes are maximum sparsity \glspl{OCT} and minimum sparsity \glspl{OCT-H} but we can specify anything in between.


\paragraph{Learning.}
With the latest computational advances in~\gls{MIO}, \cite{bd_book} were able to exactly formulate the tree training as a \gls{MIO} problem and solve it in a reasonable amount of time for problem sizes arising in real-world applications.

The \gls{OCT} cost function is a tradeoff between the misclassification error at each leaf and the tree complexity
\begin{equation*}
    \loss_{OCT} = \sum_{t \in \tree_L} L_t + \alpha \sum_{t \in \tree_B} \|a_t\|_1,
\end{equation*}
where the $L_t$ is the misclassification error at node $t$ and the second term represents the complexity of the tree measured as the sum of the norm of the hyperplane coefficients in all the splits.
The parameter $\alpha > 0$ regulates the tradeoff.
For more details about the cost function and the constraints of the problem, we refer the reader to~\cite[\Sec 2.2, \Sec 3.1]{bd_book}.

Bertsimas and Dunn apply a local search method~\citep[\Sec 2.3]{bd_book} that manages to solve \gls{OCT} problems for realistic sizes in fraction of the time an off-the-shelf optimization solver would take.
The algorithm proposed iteratively improves and refines the current tree until a local minimum is reached.
By repeating this search from different random initialization points the authors compute several local minima and then take the best one as the resulting tree.
This heuristic showed remarkable performance both in terms of computation time and quality of the resulting tree becoming the algorithm included in the OptimalTrees.jl Julia package~\citep{bd_book}.

\section{Neural Networks}
\label{sec:neural_networks}

\Glspl{NN} have recently become one of the most prominent machine learning techniques revolutionizing fields such as speech recognition~\citep{hinton2012} and computer vision~\citep{krizhevsky2012}.
The wide range of applications of these techniques recently extended to   autonomous driving~\citep{bojarski2016} and reinforcement learning~\citep{silver2017}.
The popularity of neural networks is also due to the widely used open-source libraries learning on CPUs and GPUs coming from both academia and industry such as
 TensorFlow~\citep{abadi2015}, Caffe~\citep{jia2014} and PyTorch~\citep{paszke2017}.
We use feedforward neural networks which offer   a good tradeoff between simplicity and accuracy without resorting to more complex architectures such as convolutional or recurrent neural networks~\citep{lecun2015}.

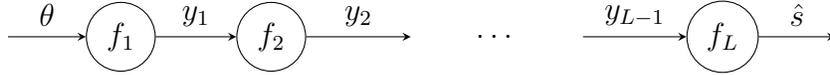
\begin{figure}
\centering
\begin{tikzpicture}[>=stealth]
\node [neuron] (f1) at (0,0) {$f_1$};
\draw [<-] (f1) -- ++(-1.5,0) node [above, midway] {$\theta$};
\node [neuron] (f2) at (2,0) {$f_2$};
\draw [->] (f1) -- (f2) node [above, midway] {$y_1$};
\node (dots1) at (4,0) {};
\draw [->] (f2) -- (dots1) node [above, midway] {$y_2$};
\node (dots) at (5,0) {$\dots$};
\node (dots2) at (6,0) {};
\node [neuron] (fL) at (8,0) {$f_L$};
\draw [->] (dots2) -- (fL) node [above, midway] {$y_{L-1}$};
\draw [->] (fL) -- ++(1.5,0) node [above, midway] {$\hat{s}$};
\end{tikzpicture}
\caption{Example Feedforward Neural Network with functions $f_i,\;i=1,\dots,L$ defined in \eqref{eq:nn_layer}.}
\label{fig:nn}
\end{figure}

\paragraph{Architecture.}
Given $L$ layers, a neural network is a composition of functions of the form
\begin{equation*}
\hat{s} = f_L(f_{L-1}(\dots f_1(\theta))),
\end{equation*}
where each function consists of
\begin{equation}
\label{eq:nn_layer}
y_{l} = f(y_{l-1}) = g(W_l y_{l-1} + b_l).
\end{equation}
The number of nodes in each layer is $n_l$ and corresponds to the dimension of the vector $y_{l} \in \reals^{n_l}$.
Layer $l=1$ is defined as the {\em input} layer and $l=L$ as the {\em output} layer. Consequently $y_{1} = \theta$ and $y_{L} = \hat{s}$.
The linear transformation in each layer is composed of an affine transformation with parameters $W_{l}  \in \reals^{n_l \times n_{l-1}}$ and $b_l \in \reals^{n_l}$.
The activation function $g: \reals^{n_l} \to \reals^{n_l}$ models nonlinearities.
We chose as activation function the~\gls{ReLU} defined as
\begin{equation*}
    g(x) = \max (x, 0),
\end{equation*}
for all the layers $l=1,\dots,L-1$.
Note that the $\max$ operator is intended element-wise.
We chose a~\gls{ReLU} because on the one hand it provides sparsity to the model since it is $0$ for the negative components of $x$ and because on the other hand it does not suffer from the vanishing gradient issues of the standard sigmoid functions~\citep{goodfellow2016}.
An example neural network can be found in Figure~\ref{fig:nn}.

For the output layer we would like the network to output not only the predicted class, but also a normalized ranking between them according to how likely they are to be the right one.
This can be achieved with a softmax activation function in the output layer defined as
\begin{equation*}
    g(x)_j = \frac{e^{x_j}}{\sum_{j=1}^{M} e^{x_j}},
\end{equation*}
where $j=1,\dots,M$ are the elements of $g(x)$.
Hence, $0 \le g(x) \le 1$ and the predicted class is $\argmax(\hat{s})$.

\paragraph{Learning.}
Before training the network, we rewrite the labels for the neural network learning using a one-hot encoding, \ie, $s_i \in \reals^M$ where $M$ is the total number of classes and all the elements of $s_i$ are 0 except the one corresponding to the class which is $1$.

We define a smooth cost function amenable to algorithms such as gradient descent, \ie, the cross-entropy loss
\begin{equation*}
    \loss_{\rm NN} =  \sum_{i=1}^{N} -s_i^\tpose \log(\hat{s}_i),
\end{equation*}
where $\log$ is intended element-wise.
The cross-entropy loss $\loss$ can also be interpreted as the distance between the predicted probability density of the labels compared to the true one.

The actual training phase consists of applying stochastic gradient descent using the derivatives of the cost function using the back-propagation rule.
This method works very well in practice and provides good out-of-sample performance with short training times.

\bibliography{refs}

\end{document}